\pdfoutput=1
\documentclass[twocolumn]{autart2}  
\usepackage{amsmath}
\usepackage{amssymb}
\usepackage{fancyhdr}
\usepackage{color}
\usepackage[latin1] {inputenc}
\usepackage{epsfig}
\usepackage{graphicx}

\def\mc{\mathcal}
\def\ee{\begin{equation}}
\def\eee{\end{equation}}
\def\bqn{\begin{eqnarray*}}
\def\eqn{\end{eqnarray*}}
\def\bnl{\begin{eqnarray}}
\def\enl{\end{eqnarray}}
\def\bma{\begin{bmatrix}}
\def\ema{\end{bmatrix}}
\def\bmx{\begin{matrix}}
\def\emx{\end{matrix}}
\def\ben{\begin{enumerate}}
\def\een{\end{enumerate}}
\def\bit{\begin{itemize}}
\def\eit{\end{itemize}}
\def\bei{\begin{itemize}}
\def\eei{\end{itemize}}
\def\bet{\begin{tabular}}
\def\eet{\end{tabular}}
\newcommand{\vn}[1]{\left|\left|#1\right|\right|}

\def\Ad{{\rm Ad}}

\def \col {\rm{col}}

\def\av{{\rm av}}

\def \RR {\mathbb{R}}

\def\qedp{\hfill{$\blacksquare$}}
\def\qed{\hfill {$\square$}}
\def\v#1{\mbox{\boldmath $#1$}}        %** bold vector **%
\def\salt{\vskip 0.4 true cm}
\newtheorem{theorem}{Theorem}

\newtheorem{definition}{Definition}

\newtheorem{remark}{Remark}
\newtheorem{proposition}{Proposition}

\begin{document}
\begin{frontmatter}

%AUTOMATICA
\title{Stabilization of Three-Dimensional Collective Motion \thanksref{footnoteinfo}} % Title, preferably not more 
\thanks[footnoteinfo]{The work is supported in part by ONR grants N00014--02--1--0826 and N00014--04--1--0534. This paper presents research results of the Belgian Network DYSCO (Dynamical Systems, Control, and Optimization), funded by the Interuniversity Attraction Poles Programme, initiated by the Belgian State,  Science Policy Office. The scientific responsibility rests with its authors.}

\author[Princeton]{Luca Scardovi}\ead{scardovi@princeton.edu},    % Add the 
\author[Princeton]{Naomi Leonard}\ead{naomi@princeton.edu},               % e-mail address 
\author[Liege]{Rodolphe Sepulchre}\ead{r.sepulchre@ulg.ac.be}  % (ead) as shown

\address[Princeton]{ Department of Mechanical and
Aerospace Engineering, Princeton University, USA}  % Please supply                                              
\address[Liege]{Department of Electrical Engineering and Computer Science,
University of Li\`ege, Belgium}             % full addresses
%\address[Baiae]{The White House, Baiae}        % here.

\begin{keyword}                           % Five to ten keywords,  
Motion coordination,
Nonlinear systems,
Multi-agent systems,
Consensus,
Multi-vehicle formations.            % chosen from the IFAC 
\end{keyword}                             % keyword list or with the 
                                          % help of the Automatica 
                                          % keyword wizard

\begin{abstract}                          % Abstract of not more than 200 words.
This paper proposes a methodology to stabilize relative equilibria
in a model of identical, steered particles moving in
three-dimensional Euclidean space. Exploiting the Lie group
structure of the resulting dynamical system, the stabilization
problem is reduced to a consensus problem on the Lie algebra. The resulting
equilibria correspond to parallel, circular and helical formations.
We first derive the stabilizing control laws in the presence of
all-to-all communication. Providing each agent with a consensus
estimator, we then extend the results to a general setting that
allows for unidirectional and time-varying communication
topologies.
\end{abstract}

\end{frontmatter}

\section{Introduction}
The problem of controlling the formation of a group of autonomous systems has received a lot of attention in recent years. This interest is principally due to the theoretical aspects that couple graph theoretic and dynamical systems concepts, and to the
vast number of applications. Applications range from sensor networks, where a group of autonomous agents has to collect
information about a process by choosing maximally informative samples \cite{LePaLeFrSeDa,bullo}, to formation control of
autonomous vehicles (e.g. unmanned aerial vehicles) \cite{JuKr_cdc,FaMu}. In these contexts it is important to consider the case
where the ambient space is the three-dimensional Euclidean space.

In the present paper we consider a model of identical particles, each with steering control, moving at unit speed in three-dimensional Euclidean space.  We address the problem of designing feedback control laws to stabilize \emph{relative equilibria} in the presence of limited communication among the agents. These equilibria are characterized by motion patterns  where the relative orientations and relative positions among the particles are constant \cite{JuKr_cdc}. The equilibria correspond to motion of all particles either 1) along parallel lines in the same direction, 2) around circles with common axis of rotation or 3) on helices with common pitch and common axis of rotation.  Therefore, our stabilization problem is a consensus problem where particles need to come to consensus on the direction, axis and pitch of their collective motion. These motion patterns are motivated by applications to vehicle groups, e.g., they provide natural and useful possibilities for collecting rich data in three-dimensional environments.  Motion patterns studied in the present paper are also motivated by the collective motion of certain animal groups \cite{CoKrJaRuFr}.

As described by Justh and Krishnaprasad \cite{JuKr_cdc}, the model for a steered, unit-speed particle can be described as a control system on the Lie group of rigid motions, $SE(3)$.   The control lives in a subspace of the Lie algebra $\frak{se}(3)$ and  provides a gyroscopic force that changes the particle's orientation (direction of motion).   Accordingly, a group of $N$ steered, unit-speed particles can be modeled as a control system on the direct product of $N$ copies of $SE(3)$.  We choose feedback control laws that depend only on relative positions and relative orientations of particles; therefore, the control preserves the $SE(3)$ symmetry of the formation. An important consequence is that no external reference is required. 

Geometry plays a central role in the investigation of the present paper and the roots of the geometric approach can be traced back to the influential work of Roger Brockett in the area of geometric control \cite{Brockett76}.  Of particular importance here, is the study of control systems on Lie groups that was formalized in Brockett's seminal work in the 1970's \cite{Brockett72,Brockett73b,Brockett73a}.  Brockett showed that system-theoretic questions, such as controllability, observability and realization theory, for a control system on a Lie group can be reduced to questions on the corresponding Lie algebra.  This work has had and continues to have enormous influence, with applications ranging from switched electrical networks \cite{Wood} to nonholonomic systems \cite{BlochBook} to control of quantum mechanical systems \cite{Navin}.

In the present paper, the geometric approach and central thesis for control systems on Lie groups are used to reduce the coordination problem on the Lie group to a consensus problem on the corresponding Lie algebra.   In particular, stabilizing particle group dynamics on $SE(3)$ is reduced to solving a consensus problem on the space of twists, $\frak{se}(3)$.

As a first step we derive stabilizing control laws
in the presence of all-to-all communication among the agents (i.e.
when each agent can communicate with all other agents at each time
instant). All-to-all communication is an assumption that is often
unrealistic in multi-agent systems. In particular, in a network of
moving agents, some of the existing communication links can fail
and new links can appear when agents leave and enter an effective
range of detection of other agents. To extend the all-to-all feedback design to
the situation of limited communication, we use the approach
recently proposed in \cite{ScSe_cdc06,scardovi:2007}, see also
\cite{FrYaLy} and \cite{Ol} for related work.

This approach suggests to replace the average quantities, often
required in a collective optimization algorithm, by local
variables obeying a consensus dynamics constrained to the
communication topology. The idea has been successfully applied to
the problem of synchronization and balancing in phase models in
the limited communication case \cite{scardovi:2007} and to the design of
planar collective motions \cite{SePaLe_limited}.

The approach leads to \emph{dynamic} control laws that include a
consensus variable that is passed to communicating
particles. The additional exchange of information is rewarded by
an increased robustness with respect to communication failures and
therefore is applicable to limited  and time-varying communication
scenarios.

On the basis of these results we design control laws that globally
stabilize collective motion patterns under mild assumptions on the
communication topology.
 
The present paper generalizes, to three-dimensional space, earlier
work in the plane \cite{SePaLe,SePaLe_limited}. Previous results in $SE(3)$ have been presented in \cite{JuKr_cdc} and in \cite{ecc,ScLeSe}. Similar approaches, applied to rigid body attitude synchronization, have been presented in \cite{SaSeLe,IgFuSp}. 
%{***\tt ADD CITATION
%SPONG***}

The rest of the paper is organized as follows. In Section
\ref{sec:model} we define the model for a group of
steered particles moving in three-dimensional Euclidean space with
unitary speed. In Section \ref{sec:all-to-all1} we review some
concepts from the theory of screws and we present a general
methodology to stabilize relative equilibria on $SE(3)$. In
Section \ref{sec:all2all} we derive control laws that stabilize
relative equilibria in the presence of all-to-all communication. 
In Section \ref{sec:consensus} we summarize some graph theoretic
notions and some results on the consensus problem in
Euclidean space. In Section \ref{sec:shape}, we design
dynamic control laws that stabilize relative equilibria in
the presence of limited communication. Finally, in Section \ref{sec:applications}, a brief discussion about possible applications in underwater robotics is presented.

For the reader's convenience the proofs of the theorems are reported in the appendix.    
%we consider the case of a limited communication topology. In this
%context we derive control laws that globally stabilize the relative equilibria of the model, by equipping each agent with a
%consensus estimator.

\section{A model of steered particles in $SE(3)$}\label{sec:model}
We consider a model of $N$ identical particles (with
unitary mass) moving in three-dimensional Euclidean space at unit
speed:
\begin{equation}\label{model_0}
\begin{array}{r c l}
\dot{\v r}_k&=&\v x_k\\
\dot{\v x}_k &=& \v u^a_k \times \v x_k, \quad k=1,2,\ldots,N,
\end{array}
\end{equation}
where $\v r_k \in \RR^3$ denotes the position of particle $k$, $\v x_k$ is the unit-norm velocity vector and $\v u^a_k \in \RR^3$
is a control vector. Model (\ref{model_0}) characterizes particle dynamics with forcing only in the directions normal to
velocity, i.e. $\ddot{\v r}_k = \v u^a_k \times \dot{\v r}_k$. An alternative to (\ref{model_0}) is to provide each particle
with an orthonormal frame and to write the system dynamics in a \emph{curve framing} setting \cite{JuKr_cdc}:
\begin{equation}\label{model}
\begin{array}{r c l}
\dot{\v r}_k&=&\v x_k\\
\dot{\v x}_k&=&\v y_k q_k + \v z_k h_k\\
\dot{\v y}_k&=&-\v x_k q_k + \v z_k w_k \\
\dot{\v z}_k&=&-\v x_k h_k - \v y_k w_k, \quad \quad k=1,\ldots,N,\\
\end{array}
\end{equation}
where $(\v x_k,\v y_k,\v z_k)$ is a right handed orthonormal frame associated to particle $k$ (in particular $\v x_k \in S_2$ is
the (unit) velocity vector). The scalars $q_k$, $h_k$ represent the curvature controls of the $k$th particle. The scalar $w_k$
adds a further degree of freedom allowing rotations about the axis $\v x_k$. In vector notation we define
\begin{equation}\label{u}
\v u_k =\left[
\begin{matrix} w_k\\-h_k\\q_k
\end{matrix}\right].
\end{equation}
%It is worth noting that in (\ref{model}), there is freedom in the choice of initial conditions $\v y_k(0)$ and $\v z_k(0)$.
%In the following we will use a bold variable without index to denote $3N$-vectors, e.g. $\v x = (\v x_1^T, \ldots, \v x_N^T)^T$.
%The notation $\bar{\v x}$ will indicate an average vector, i.e.
%\[
%\bar{\v x} = \frac{1}{N}\sum_{k=1}^N \v x_k.
%\]
The advantage of using model (\ref{model}) instead of model (\ref{model_0}) relies on its group structure. Model (\ref{model})
indeed defines a control system on the Lie group $SE(3)$ and the dynamics (\ref{model}) can be expressed in terms of the group
variables $g_k \in SE(3)$:
\begin{equation}\label{model_lie}
\dot g_k = g_k \hat{\v \xi}_k, \quad \quad k=1,\ldots,N,
\end{equation}
where $\hat{\v \xi}_k \in \frak{se}(3)$ is an element of the Lie algebra
of $SE(3)$, the tangent space to $SE(3)$ at the identity. From (\ref{model}) we obtain
\[
\begin{array}{c c}
g_k = \left[\begin{matrix} R_k& \v r_k\\\v 0& 1 \end{matrix}\right],& R_k=\left[\v x_k, \v y_k, \v z_k  \right] \in SO(3),
\end{array}
\]
and
\begin{equation}\label{lie1}
\hat{\v \xi}_k = \left[
\begin{matrix} \hat{\v u}_k& \v e_1\\\v 0& 0
\end{matrix}\right],
\end{equation}
where
%\[
%R_k=\left[\v x_k, \v y_k, \v z_k  \right] \in SO(3)
%\]
%and
\[
\hat{\v u}_k =\left[
\begin{matrix} 0& -q_k & -h_k\\q_k&0&-w_k\\h_k&w_k&0
\end{matrix}\right]
\]
is a skew-symmetric matrix that represents an element of $\frak so(3)$, the Lie algebra of $SO(3)$. We denote by $(\v e_1, \v e_2,\v e_3)$
the standard orthonormal basis for $\RR^3$.

When only the orientations of the particles are taken into account, the reduced dynamics of (\ref{model_lie}) are
\begin{equation}\label{model3}
\dot{R}_k=R_k\ \hat{\v u}_k, \quad \quad k=1,\ldots,N
\end{equation}
and the system evolves on the Lie group $SO(3)$.

It is worth noting that the following relation exists between the control vector $\v u^a_k$ in (\ref{model_0}) and the vector
$\v u_k$ in (\ref{u}):
\begin{equation}\label{curvat}
\v u^a_k = R_k \v u_k.
\end{equation}
Therefore $\v u^a_k$ can be interpreted as the control vector $\v u_k$ expressed in the \emph{spatial} reference frame\footnote{We adopt the word spatial to mean ``relative to a fixed (inertial) coordinate frame''.}.

If the curvature controls in model (\ref{model}) are feedback functions of \emph{shape} quantities (i.e. relative frame orientations and relative positions), the closed-loop vector field is invariant under an action of the symmetry group $SE(3)$.
The resulting closed-loop dynamics evolve in a quotient manifold called \emph{shape space} and the equilibria of the reduced dynamics are called \emph{relative equilibria}. To formally introduce the shape variable associated to two particles $k$ and $j$ we define
\[
g_{kj}\triangleq g_k^{-1} g_j
\]  
 which, in the case of dynamics evolving on $SE(3)$, particularize to
 \[ 
 g_{kj} = 
 \left[
 \begin{matrix}
 R_{kj} & r_{jk}^{k}\\
 0& 1
 \end{matrix}
\right]
\]
where $R_{kj}\triangleq R_{k}^{T}R_{j}$ and $r_{jk}^{k}\triangleq  R_k^T\left(r_j-r_k\right)$. 
As pointed out previously, our control laws will be restricted to depend on shape variables only. 
Therefore the (static and dynamic) control laws will assume the form   
\[
u_{k}= \eta^{s}_{k}({\mathcal R}_{kj},{\mathcal D}_{jk}^{k}), 
\] 
and 
\[
\begin{array}{rcl}
u_{k}&=& \eta^{d}_{k}({\mathcal R}_{kj},{\mathcal D}_{jk}^{k},\gamma_{k})\\
\dot \gamma_{k}&=&\rho_{k}({\mathcal R}_{kj},{\mathcal D}_{jk}^{k},\gamma_{k}),
\end{array}
 \]
respectively, where ${\mathcal R}_{kj} = \{R_{kj},\, j=1,2,\ldots,N\}$,\, ${\mathcal D}_{jk}^{k}= \{r_{kj}^{k},\, j=1,2,\ldots,N\}$ and $\gamma_{k}$ are additional consensus variables.
As we will see in the following, dynamic control laws will be used several times in the paper.
In particular it will turn out that (in general) to stabilize relative equilibria in a decentralized framework a static control law is not sufficient. Furthermore, as pointed out in earlier works \cite{scardovi:2007,ScSe_cdc06,SePaLe_limited}, dynamic control laws are required when a limited communication setting is taken into account (see Section \ref{sec:shape}).
         
\begin{figure}[t]
      \centering
      \begin{tabular}{l}
  \hspace{-0.55 cm}  \includegraphics[scale =0.55]{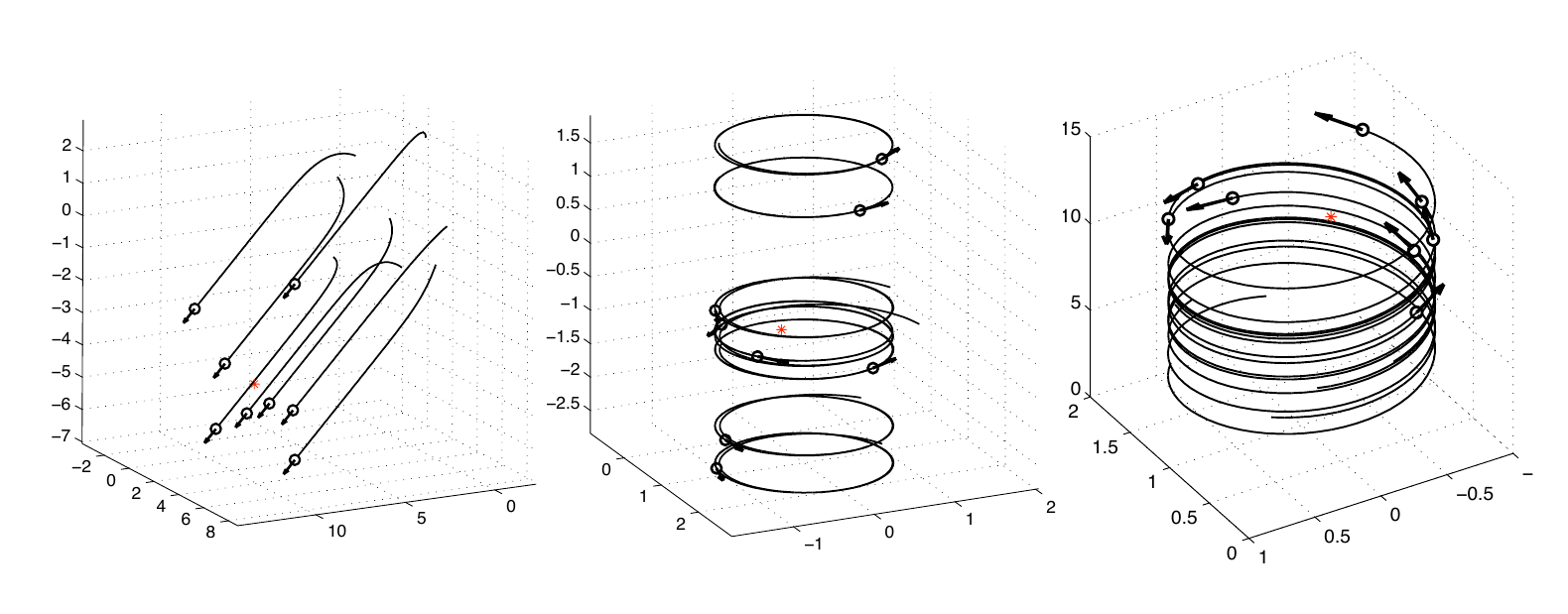}\\
 \hspace{1 cm} (a)   \hspace{2.2 cm} (b)  \hspace{2 cm} (c)
 %\hspace{-1.8 cm}\includegraphics[scale=0.4]{circ.eps}&\hspace{-2.8 cm} \includegraphics[scale=0.4]{helix.eps}
      \end{tabular}
      \caption{The three types of relative equilibria: (a) parallel, (b) circular and (c) helical.
      }
      \label{fig:formations}
   \end{figure}

Relative equilibria of the model (\ref{model}) have been characterized in \cite{JuKr_cdc}. 
The equilibria, depicted in Figure \ref{fig:formations}, are of three types:\\
i) Parallel motion: all particles move in the same direction with arbitrary relative positions;\\
ii) Circular motion: all particles draw circles with the same radius, in planes orthogonal to the same axis of rotation;\\
iii) Helical motion: all particles draw circular helices with the same radius, pitch, axis and axial direction of motion. 

In the following section we will show how to characterize the relative equilibria by using screw theory. This approach will be particularly useful in Section \ref{sec:all2all} when 
the problem of stabilizing the relative equilibria will be addressed.
%In \cite{ecc} the first two relative equilibria have been stabilized by means of \emph{static} control laws in presence of
%all-to-all communication. In presence of limited communication \emph{dynamic} control laws  (equipping each particle with a
%consensus estimator) provide the desired convergence results. Unfortunately, the stabilization of the helical equilibria by
%means of static control appears to be more challenging.
%In this paper we follow a different approach. We identify a potential function (based on a common reference vector) whose global
%minima characterize relative equilibria of (\ref{model}). Unfortunately the obtained control laws will not depend only on the
%\emph{shape} of the formation and therefore the resulting closed loop dynamics will not be invariant under an action of the
%symmetry group $SE(3)$. To recover this symmetry, in Section (\ref{shape}) we propose dynamic control laws that are based on a
%consensus dynamics.

%\begin{figure}[t]
%      \centering
%      \begin{tabular}{c c c}
%\hspace{-1 cm}      \includegraphics[scale=0.35]{parallel.eps}&\hspace{-3 cm}    \includegraphics[scale=0.35]{circ.eps}& \hspace{-3 cm}    \includegraphics[scale=0.35]{helix.eps}\\
%\hspace{-1.5 cm}      (a)& \hspace{-3 cm} (b)& \hspace{-3 cm}(c)
% %\hspace{-1.8 cm}\includegraphics[scale=0.4]{circ.eps}&\hspace{-2.8 cm} \includegraphics[scale=0.4]{helix.eps}
%      \end{tabular}
%      \caption{The three types of relative equilibria: (a) parallel, (b) circular and (c) helical.
%      }
%      \label{fig:formations}
%   \end{figure}
   
\section{Stabilization of relative equilibria as a consensus problem}\label{sec:all-to-all1}
In terms of screw theory \cite{MuZeSh}, an element of $\frak{se}(3)$ is called %$\hat{\xi}_k$ is called
a \emph{twist}. The motion produced by a constant twist is called a \emph{screw motion}. The operator denoted by $\vee$ extracts
the $6$-dimensional vector which parameterizes a twist: (\ref{lie1}) yields
\[
\v \xi_k =
 \left[
\begin{matrix} \hat{\v u}_k& \v e_1\\\v 0& 0
\end{matrix}\right]^{\vee} = \left[
\begin{array}{c}
\v e_1\\
\v u_k
\end{array}
 \right].
\]
The inverse operator, $\wedge$, expresses the twist in homogeneous coordinates starting from a vector form: (\ref{lie1}) yields
\[
\hat{\v \xi}_k
 = \left[
\begin{array}{c}
\v e_1\\
\v u_k
\end{array}
 \right]^{\wedge} =  \left[
\begin{matrix} \hat{\v u}_k& \v e_1\\\v 0& 0
\end{matrix}\right].
\]
A constant twist $\v \xi_0 = [\v v_0^T,\v \omega_0^T]^T \in \RR^6$ defines the
screw motion $g(0)e^{\hat{\v \xi_0}t}$ on $SE(3)$ \cite{MuZeSh}, where $g(0)$ denotes the initial condition. When
$\v \omega_0 \neq 0$ this motion yields a final configuration that corresponds to a rotation by the
amount $\theta=\vn{\v \omega_0}$ about an axis $l$, followed by
translation by an amount $p_0\vn{\v \omega_0}$ parallel to the axis
$l$. When $\v \omega_0 = 0$ the corresponding screw motion
consists of a pure translation along the axis $\lambda \v v_0$ of the screw by a
distance $M_0 = \vn{\v v_0}$.
%A screw motion is defined as follows \cite{MuZeSh} \salt
%\begin{definition}
%A screw consists of an axis $\v l$, a pitch $p$ and a magnitude $M$. A screw motion represents rotation by an amount $\theta=M$
%about the axis followed by translation by an amount $h\,\theta$ parallel to the axis $\v l$. If $p=\infty$ then the
%corresponding screw motion consists of a pure translation along the axis of the screw by a distance $M$.
%\end{definition}\salt
%Let
%\[
%\v \xi = [\v v^T,\v \omega^T]^T \in \RR^6
%\]
%be a twist.
The relations among the screw $(l_0, p_0, M_0)$ and twist $\v \xi$ are the
following \cite{MuZeSh}:\\
%\[
%p = \frac{\v \omega^T \v v}{\vn{\v \omega}^2}
%\]
\[
\begin{array}{rcl}
p_0&=&\left\{
\begin{array}{l r}
\frac{\v \omega_0^T \v v_0}{\vn{\v \omega_0}^2},& \;\;\quad \quad \mbox{if}\; \v \omega_0 \neq 0\\
\infty, & \mbox{if}\; \v \omega_0 = 0
\end{array}
\right.\\
l_0&=&\left\{
\begin{array}{l r}
\frac{\v \omega_0 \times \v v_0}{\vn{\v \omega_0}^2} + \lambda \v
\omega_0 ,& \mbox{if}\; \v \omega_0 \neq 0\\
0 + \lambda \v v_0, & \mbox{if}\; \v \omega_0 = 0
\end{array}
\right.\\
 M_0&=&\left\{
\begin{array}{l r}
\vn{\v \omega_0} ,&\quad \quad \quad \mbox{if}\; \v \omega_0 \neq 0\\
\vn{\v v_0}, &\quad \quad \quad \mbox{if}\; \v \omega_0 = 0
\end{array}
\right.
\end{array}
\]
where $\lambda \in \RR$. 

In the context of model (\ref{model_lie}), the twist (in body coordinates) is given by $\v \xi_k=[\v
e_1^T, \v u_k^T]^T$.
%\[
%\v \xi_k = \left[
%\begin{array}{c}
%\v e_1\\
%\v u_k
%\end{array}
% \right].
%\]
To map $\v \xi_k$ into a \emph{spatial} reference frame, one uses the \emph{adjoint transformation} associated with $g_k$
\[
\Ad_{g_k} = \left[
\begin{matrix}
R_k & \hat{\v r}_k R_k\\
0& R_k
\end{matrix}
\right],
\]
which yields
\begin{equation}\label{twist_glob}
%\begin{array}{rcll}
\v \xi_k^a \triangleq \Ad_{g_k} \v \xi_k = \left[\begin{array}{c}
\v x_k + \v r_k \times R_k \v u_k\\
R_k \v u_k
\end{array}\right]\hspace{-0.1 cm}=\hspace{-0.1 cm}\left[\begin{array}{c}
\v x_k + \v r_k \times {\v u}_k^a\\
{\v u}_k^a
\end{array}\right].
%\end{array}
\end{equation}
To give a geometric interpretation to (\ref{twist_glob}) we
compute the relative screw coordinates (expressed in the spatial
frame) and we obtain an (instantaneous) pitch
\[
p_k = \frac{\v e_1^T \v u_k}{\vn{\v u_k}^2},
\]
an (instantaneous) axis
\[
l_k^a=\left\{
\begin{array}{l r}
{\v u}^a_k \times \frac{\v x_k + \v r_k \times {\v u}^a_k}{\vn{\v u_k}^2} + \lambda \v
u^a_k ,& \mbox{if}\; \v u_k \neq 0\\
0 + \lambda \v x_k, & \mbox{if}\; \v u_k = 0
\end{array}
\right.,
\]
and (instantaneous) magnitude
\[
M_k=\left\{
\begin{array}{l l}
\vn{\v u_k} ,& \mbox{if}\; \v u_k \neq 0\\
1, & \mbox{if}\; \v u_k = 0.
\end{array}
\right.
\]
Therefore, constant control vectors $\v u_k$, $k=1,2,\ldots,N$, define screw motions (corresponding to helical, circular or straight motions).
%the following types of motion\\
%i) if $\vn{\v u_k} = 0$: straight motion along the line $\lambda \v x_k$, $\lambda \in \RR$\\
%ii) if $\vn{\v u_k} \neq 0$ and $w_k = 0$: circular motion about the axis ${\v u}^a_k \times \frac{\v x_k + \v r_k
%\times {\v u}^a_k}{\vn{\v u_k}^2} + \lambda \v u^a_k$, $\lambda \in \RR$, with radius $\vn{\v u_k}^{-1}$.\\
%iii) if $\vn{\v u_k} \neq 0$ and $w_k \neq 0$: helical motion about the axis ${\v u}^a_k \times \frac{\v x_k + \v r_k \times {\v
%u}^a_k}{\vn{\v u_k}^2} + \lambda \v u^a_k$, $\lambda \in \RR$, with pitch $\displaystyle \frac{w_k}{\vn{\v u_k}^2}$ and radius
%$\sqrt{1-w_k^2\vn{\v u_k}^2}/\vn{\v u_k}$.\salt

%\section{Relative Equilibria and Consensus on Twists}\label{all-to-all1}
Now we are ready to geometrically characterize the relative equilibria of (\ref{model_lie}).
Consider two particles and their respective group variables $g_k$ and $g_j$. The dynamics for $g_{kj}= g_k^{-1} g_j$ (the shape
variable) are given (see \cite{JuKr_cdc}) by
\begin{equation}\label{mod_shape}
\begin{array}{rcl}
\dot{g}_{kj} &=& -g_k ^{-1} \dot g_k  g_k ^{-1} g_j+g_k ^{-1} g_j\hat{\v \xi}_j\\
&=&-\hat{\v \xi_k} {g}_{kj} + {g}_{kj}\hat{\v \xi}_j\\
&=&g_{kj}(\hat{\v \xi}_j - \widehat{\Ad_{g_{kj}^{-1}} \v \xi_k}).
\end{array}
\end{equation}
%The equilibria of (\ref{mod_shape}) are relative equilibria of (\ref{model_lie}) (for two particles).
Equation (\ref{mod_shape}) implies that a relative equilibrium of
(\ref{model_lie}) is reached when the twists (expressed into a spacial reference frame) are equal for all the particles, i.e. $\v
\xi_k^a = \v \xi^a_0$ for $k=1,\ldots,N$, $\v \xi^a_0$ arbitrary.
To see it, it is sufficient to  equate the last term in
(\ref{mod_shape}) with zero  and to apply the
adjoint transformation $\Ad_{g_j}$ obtaining
\begin{equation}\label{mod_shape2}
\begin{array}{rcl}
\Ad_{g_j} {\v \xi}_j - \Ad_{g_j} \Ad_{g_{kj}^{-1}} {\v
\xi}_k= {\v \xi}^a_j - {\v \xi}^a_k = 0.
\end{array}
\end{equation}

Since the screw coordinates associated to the common value $\v \xi^a_0$
provide a geometrical description of the motion, the relative equilibria are characterized by a pitch, an axis and a magnitude uniquely determined by $\v \xi^a_0$.
We summarize the above discussion in the following Proposition. Let ${\bf 1}_{N} = (1,\ldots,1)^{T}\in \RR^{N}$.
\begin{proposition}\label{prop}
The following statements are equivalent:\\
i) System (\ref{model}) is at a relative equilibrium.\\
ii) The twists $\v \xi_k^a$ defined by (\ref{twist_glob}) are equal for $k=1,2,\ldots,N$, i.e. the following algebraic condition
is satisfied
\[
\tilde{\Pi} \v \xi^a = 0,
\]
where $\tilde{\Pi} = (I_{N}-\frac{1}{N}\v 1_{N} \v 1_{N}^T) \otimes I_6$ and ${\v \xi}^a = {\col}({\v {\xi}}^a_{1},\ldots,\v{\xi}_N^a)$.\qed 
\end{proposition} 

Proposition \ref{prop} reduces the problem of stabilizing a relative equilibrium on $SE(3)$ to a consensus problem on
twists.  

In the rest of the paper, we denote by $\Sigma$ the set of solutions of (\ref{model}) with consensus on the rotation vector, i.e. $\v  \omega_{k}^{a} = \v \omega_j^a,\, k,j=1,2,\ldots,N$ :
\[
\begin{array}{rcl}
\Sigma&\triangleq& \{g_{k} \in SE(3),\,k=1,\ldots,N:\; g_{k} = g_{k}(0)e^{\hat{\v \xi}_{k} t},\\
&& \v \xi_{k} ={\col}(\v e_{1},R_{k}^{T} \v \omega^a_k),\, \v  \omega_{k}^{a} = \v \omega_j^a, j=1,2,\ldots,N,\\ 
&&g_{k}(0) \in SE(3)\} 
\end{array}
\]
and we denote by $E$ the subset of $\Sigma$ corresponding to relative equilibria. By Prop. \ref{prop}, this set is characterized as 
\[
\begin{array}{rcl}
E& \triangleq& \{g_{k} \in SE(3),\,k=1,\ldots,N: g \in \Sigma,\\
 &&\v v_k^a = \v v_j^a,\, j=1,2,\ldots,N \}.   
\end{array}
\]

Likewise we will denote by $\Sigma(\v \omega_0)$ the subset of $\Sigma$ where $\v \omega^{a}_k = \v \omega_0,\, k=1,2,\ldots,N$, for some fixed $\v \omega_0 \in \RR^3$ and by $E(\v \omega_0)$ the subset of $E$ with a fixed rotation vector $\v \omega_0$. 
%\[
%\begin{array}{rcl}
%E(\v \omega_0,\v v_0)& \triangleq& \{g_{k} \in SE(3),\,k=1,\ldots,N: g \in \Sigma(\v \omega_0),\\
% &&\v v_k^a =\v v_0\}.   
%\end{array}
%\]
%where $\v v_0 \in \RR^3$, the subset of $E(\v \omega_0)$ that defines a screw specified by the twist $\xi_0^a = {\col}(\v v_0, \v \omega_0)$.     
\begin{remark}
The discussion above particularizes to $SE(2)$. Consider the (planar) model 
\begin{equation}\label{modelse2}
\begin{array}{r c l}
\dot{\v r}_k&=&\v x_k\\
\dot{\v x}_k&=&u_k \v y_k  \\
\dot{\v y}_k&=&-u_k \v x_k. \\
\end{array}
\end{equation}
for $k=1,\ldots, N$.
In the Lie group $SE(2)$, we obtain
\[
\begin{array}{c c}
g_k = \left[\begin{matrix} R_k& \v r_k\\ 0& 1
\end{matrix}\right], &
 \hat{\v \xi}_k = \left[
\begin{matrix}  \hat{\v u}_k& \v e_1\\ 0& 0
\end{matrix}\right]
\end{array}
\]
for  $k=1,\ldots, N$, where
\[
R_k=\left[\v x_k, \v y_k\right] \in SO(2),
\]
\[
\hat{\v u}_k = \left[\begin{matrix}  0& -u_k\\ u_k& 0
\end{matrix}\right]=Ju_k,\quad \quad J = \left[\begin{matrix} 0 &-1\\1&0 \end{matrix}\right]
\]
and $\v e_1 = [1,0]^T$. In this case the twist is $\v \xi_k = [\v e_1^T,u_k]^T \in \RR^3$. By mapping the twist coordinates to a
spatial frame we obtain
\begin{equation} \label{twse2}
\v \xi^a_k = \left[
\begin{array}{c}  \v x_k - u_k J \v r_k\\
u_k
\end{array}
\right], \quad \quad  k=1,\ldots, N.
\end{equation}
When $u_k,  k=1,\ldots, N,$ are constant, only two types of motion are possible for (\ref{modelse2}), straight motion ($u_k=0$) and circular
motion ($u_k=\omega_0$). When (\ref{twse2}) are equal and constant for all the particles the resulting motion is characterized
by a parallel formation ($u_k=0$) and a circular formation about the same point ($u_k\neq 0$ and constant). Stabilizing control laws are
derived in \cite{SePaLe,SePaLe_limited}.
\end{remark}
%The foregoing result provides, in a very natural way, a design
%methodology to stabilize the relative equilibria of (\ref{model}).
%In fact, the problem of stabilizing a relative equilibrium can be
%reduced to a consensus problem on the twists.
%It is worth noting that a relative equilibrium is characterized by a configuration in which the twists expressed in a common
%reference frame (or equivalently the screw parameters) are equal (and constant) for all the particles. Therefore, the problem of
%stabilizing the relative equilibria can be reduced to a consensus problem on the twists.
\section{Stabilization of relative equilibria in the presence
of all-to-all communication}\label{sec:all2all}
From (\ref{twist_glob}), when a particle $k$ applies the constant control $\v u_k = \v \omega_k$, the (constant) twist expressed
in the spatial reference frame is
\begin{equation}\label{twi}
\v \xi^a_k=
\left[\begin{array}{c}
\v x_k + \v r_k \times R_k \v \omega_k\\
R_k \v \omega_k
\end{array}\right] =\left[\begin{array}{c}
\v v_k^a\\
\v \omega^a_k
\end{array}\right].
\end{equation}
Motivated by Proposition \ref{prop} a natural candidate Lyapunov function is
\begin{equation}\label{V}
V(\v \xi^a) =\frac{1}{2}\vn{\tilde{\Pi} \v \xi^a}^2=\frac{1}{2}\sum_{k=1}^N \vn{\v \xi_k^a - \v \xi^a_{\av}}^2
\end{equation}
where the subscript ``$\av$'' is used to denote average quantities, i.e.
\[
{\v \xi_{\av}^a} = \frac{1}{N} \sum_{k=1}^N \v \xi_k^a.
\]
This is the approach pursued in \cite{SePaLe} for collective motion in $SE(2)$.

Unfortunately, from (\ref{twi}), it is evident that the first
component $\v v_k^a$ is not linear in the state variables. As a
consequence $\v v^a_{\av} \neq \v x_{\av} +  \v r_{\av}\times \v
\omega^a_{\av}$ and the approach followed in \cite{SePaLe} does not
yield shape control laws. To understand how to overcome this obstacle we first
stabilize the motion about an axis of rotation with direction that is
fixed. In Section \ref{sec:consome} we relax the design by
replacing, in the control laws, the fixed direction of the axis of
rotation with (local) consensus variables, thereby obtaining
stabilizing shape control laws. A simplification occurs when the
desired relative equilibrium corresponds to parallel formations.
For this relative equilibrium the twists reduce to the velocity
vectors and therefore a simplified consensus problem may be
addressed. In the next section we address this simpler case, while
the general case is addressed in Section \ref{sec:fix} and in
Section \ref{sec:consome}. 

\subsection{Stabilization of parallel formations}\label{sec:parallel}
First observe that when the particles follow straight trajectories
(\ref{twi}) reduces to
\[
\v \xi^a_k= \left[\begin{array}{c}
\v x_k \\
\v 0
\end{array}\right], \quad \quad  k=1,\ldots, N,
\]
and the Lyapunov function (\ref{V}) reduces to
\begin{equation}\label{V_x}
V(\v x)=\frac{N}{2}\left(1-\vn{{\v x}_\av}^2\right).
\end{equation}
%Notice that minimizing (maximizing) (\ref{W_x}) is equivalent to
%maximizing (minimizing) the potential
%\begin{equation}\label{V_x}
%V(\v x)=\frac{N}{2}\vn{{\v x}_\av}^2.
%\end{equation}
The parameter $\vn{{\v x}_{\av}}$ is a measure
of synchrony of the velocity vectors $\v x_k,\, k=1,2,\ldots,N$.
In the model (\ref{model}), $\vn{{\v x}_{\av}}$ is maximal when
the velocity vectors are all aligned (synchronization) leading to
parallel formations. It is minimal when the velocities balance to
result in a vanishing centroid, leading to collective motion
around a fixed center of mass. Synchronization (balancing) is therefore
achieved by minimizing (maximizing) the potential (\ref{V_x}). The
time derivative of (\ref{V_x}) along the solutions of
(\ref{model3}) is
\begin{equation}\label{V_x_dot}
%\begin{array}{rcl}
\dot{V}= -\sum_{j=1}^N <{\v x}_{\av},\dot{\v x}_j> 
%&=& \sum_{j=1}^N
%<{\v x}_{\av},\v y_j>q_j + <{\v x}_{\av},\v z_j>h_j,
%\end{array}
\end{equation}
where $<\cdot,\cdot>$ denotes the scalar product.

The control law
%\begin{equation}\label{cont1}
%\begin{array}{rcl}
%q_k&=&-K<{\v x}_{\av}, \v y_k> \\
%h_k&=&-K<{\v x}_{\av}, \v z_k>, \quad \quad k=1,\ldots,N,
%\end{array}
%\end{equation}
\begin{equation}\label{cont1}
\v u_k=R_k^T(\v x_k \times \v x_{\av}), \quad \quad  k=1,\ldots, N,
\end{equation}
ensures that (\ref{V_x}) is non-increasing.\\
%In vector form (\ref{cont1}) reduces to
%\begin{equation}\label{par_vec}
%\v u_k=-K R_k^T(\v x_k \times \v x_{\av}).
%end{equation}
 The following result provides a characterization of the dynamics 
of model (\ref{model}) with the control law (\ref{cont1}).
\begin{theorem}\label{all-to-all}
Consider the model (\ref{model}) with the control law
(\ref{cont1}). The closed-loop vector field is invariant under an action of
the group $SE(3)$. Every solution exists for all $t \ge 0$ and asymptotically
converges to $\Sigma(0)$. Furthermore, the set $E(0)$ of parallel motions is asymptotically stable
in the shape space and every other positive limit set is unstable. \qed
\end{theorem}

\begin{figure}[thb]
      \centering
      \begin{tabular}{c c}
      \hspace{-0.45 cm}\includegraphics[scale=0.30]{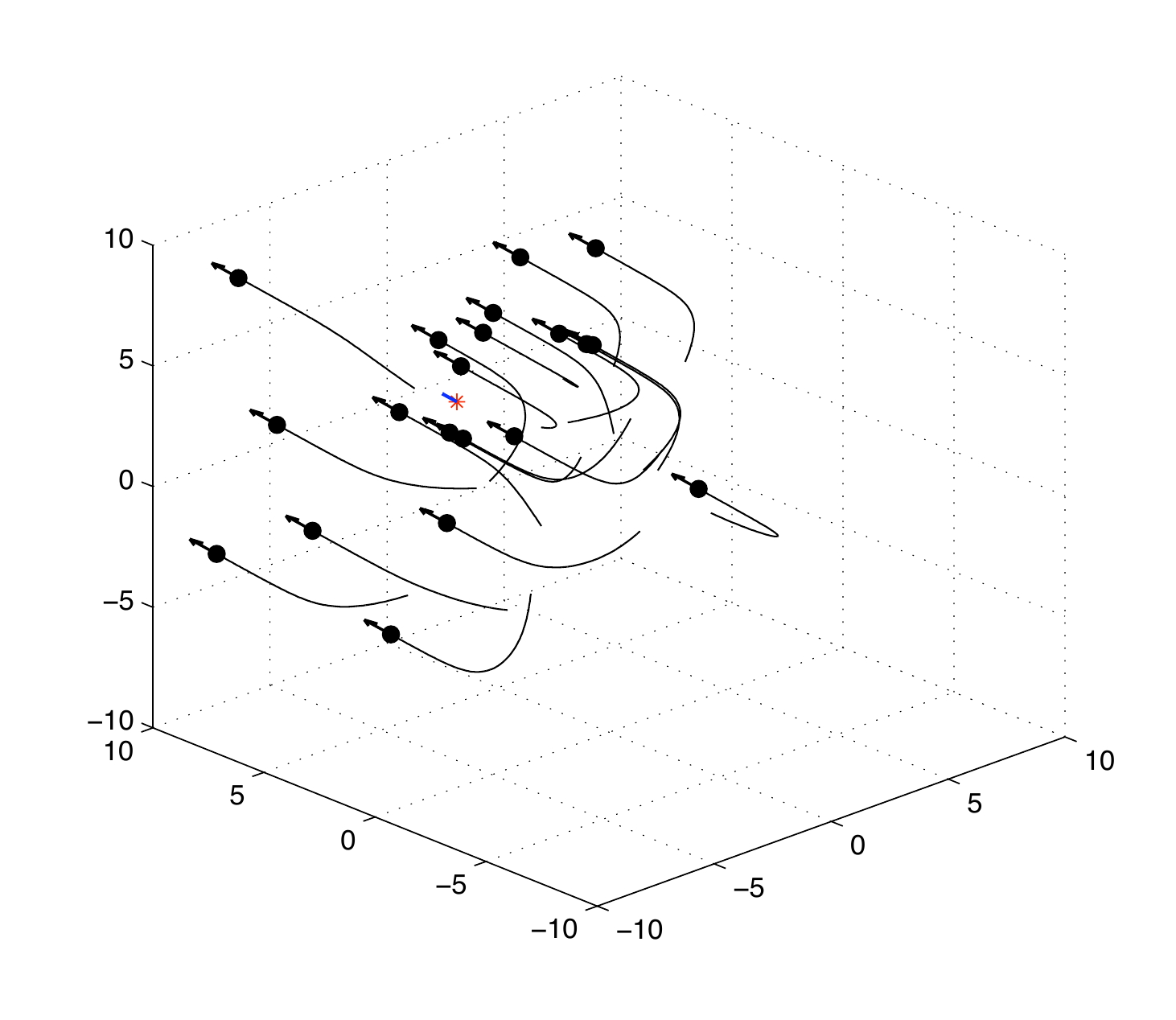}&\hspace{-0.8 cm}\includegraphics[scale=0.34]{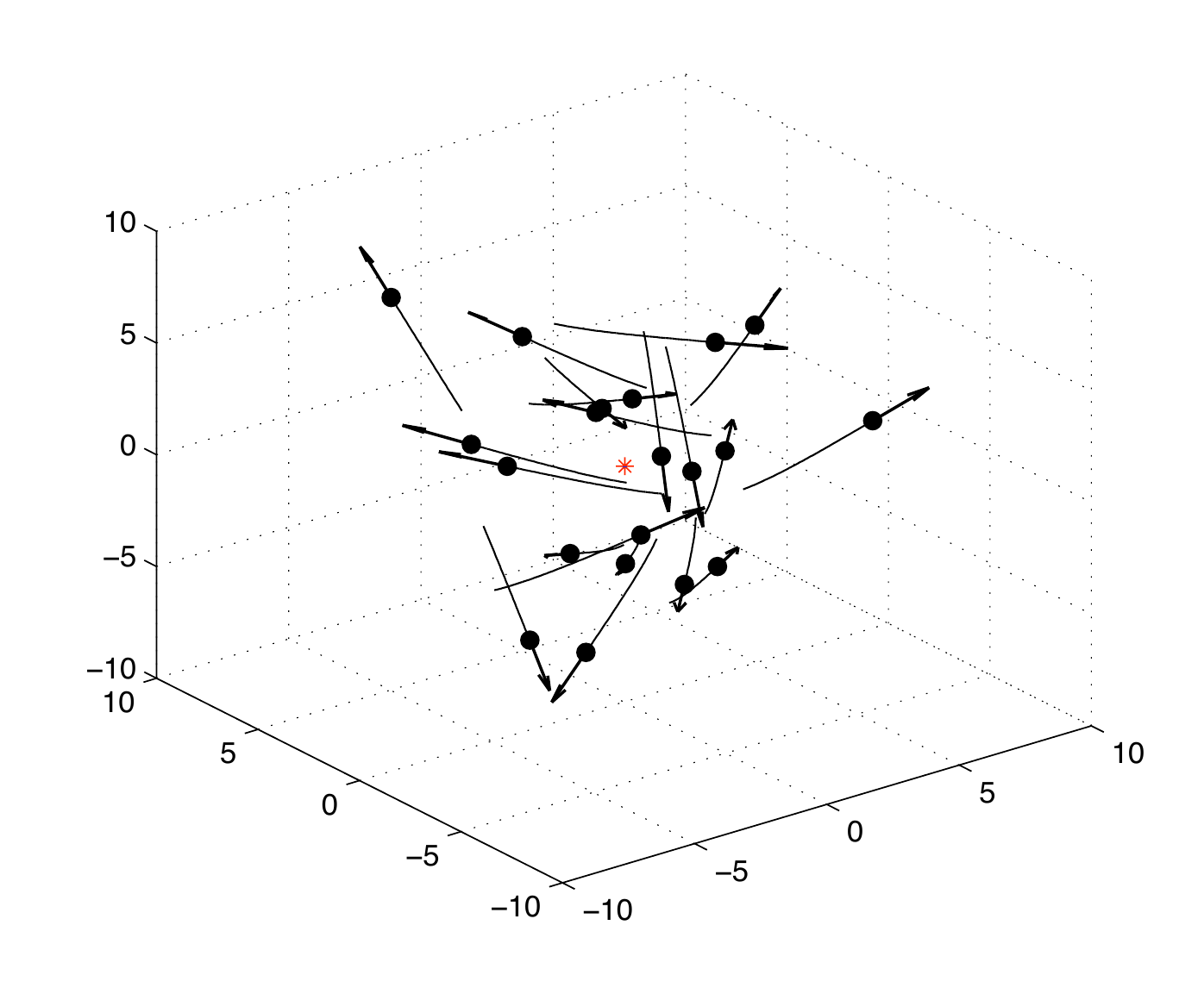}\\
   \hspace{-0.45 cm}    (a)& \hspace{-0.8 cm}(b)
      \end{tabular}
      \caption{Parallel and balanced formations.}
      \label{fig:sim1}
   \end{figure}

As a consequence of Theorem~\ref{all-to-all}, we
obtain that the control law (\ref{cont1}) stabilizes parallel
formations (see Fig.~\ref{fig:sim1}a). 
\begin{remark}
When the sign is reversed in (\ref{cont1}), only the set of balanced states (i.e. those states such that ${\v x}_{\av}$ is zero) is asymptotically stable and every other
equilibrium is unstable.
This leads to configurations where the center of mass of the particles is a
fixed point (see Fig.~\ref{fig:sim1}b). The stabilization
of the center of mass to a fixed point does not lead in general to
a relative equilibrium and therefore is not of interest in the present paper.
\end{remark}
\begin{remark}
It is worth noting that the feedback control (\ref{cont1}) does
not depend on the relative orientation of the frames but only on
the relative orientations of the velocity vectors. Therefore, each particle compares only relative velocity vectors with respect to its own reference frame, in order to implement
control law (\ref{cont1}).
\end{remark}

\subsection{Stabilization of screw relative equilibria: preliminary design}\label{sec:fix}
Let $\v \omega_0 \in \RR^3$ be a fixed
constant vector expressed in the spacial reference frame. Observe that under the constant control law $
\v u_{k} = R_k^T \v \omega_0$, a relative equilibrium is reached
when the vectors $\v v_k^a$ in (\ref{twi}) are equal for all the
particles.

Up to an additive constant the Lyapunov function (\ref{V}) becomes
\begin{equation}\label{S}
S(\v v^a, \v \omega_0) = \frac{1}{2}\sum_{k=1}^N\vn{\v v_k^a - \v v^a_{\av}}^2
\end{equation}
where $\v v_k^a = \v x_k + \v r_k \times \v \omega_0$ and $\v v^a = {\col}(\v v^a_1,\ldots,\v v^a_N)$. The time derivative is
\[
\dot{S}=\hspace{-1 mm} \sum_{k=1}^N <\v v_k^a - \v v^a_{\av},\, \dot{\v v}^a_k> =\hspace{-1 mm}\sum_{k=1}^N <\v v_k^a - \v v^a_{\av}, \dot{\v x}_k + \v x_k \times \v
\omega_0>.
\]
The control law
\begin{equation}\label{concon2}
\v u_k =  R^T_k \left (\v \omega_0 + \left[ \left(\v r_k - \v r_{\av}\right)\times \v \omega_0 - \v x_{\av}\right]\times \v
x_k\right),
\end{equation}
for  $k=1,\ldots, N,$ results in a non-increasing $S$
\begin{equation}\label{Sdot}
\dot S = - \sum_{k=1}^{N}\vn{ \Pi_{\v x_k} (\v v_k^a-\v v^a_{\av}) }^2 \le 0.
\end{equation}
where $\Pi_{{\v x_k}} =I-\displaystyle {\v x_k \v x_k^T}$ is the projection matrix on the orthogonal complement of the subspace spanned by $\v x_k$. 
Note that the $\v v_k^a$ dynamics with the control law (\ref{concon2}) are
\begin{equation}\label{s0}
 \dot{\v v}^a_k =-\Pi_{{\v x_k}}\left(
\v v_k^a -\v v^a_{\av} \right), \quad \quad  k=1,\ldots, N.
\end{equation}
%By using the definition (\ref{u}) the control law
%(\ref{concon3}) can be written with respect its components
%\begin{equation}\label{concon2}
%\begin{array}{c c l}
%q_k &=&  <\v \omega,\v z_k> - <(\v r_k -\v r_{\av}) \times \v \omega- \v x_{\av}, \v y_k>\\
%h_k &=& -<\v \omega,\v y_k> - <(\v r_k -\v r_{\av}) \times \v \omega- \v x_{\av}, \v z_k>\\
%w_{k}&=& <\v \omega,\v x_k, >.
%\end{array}
%\end{equation}
The convergence properties of the resulting closed-loop system are characterized in the following theorem:% (see Appendix B for the proof):     
\begin{theorem}\label{th:1}
Consider model (\ref{model}) with the control law (\ref{concon2}). The closed-loop vector field is invariant  under an action of
the translation group $\RR^3$. Every solution exists for all $t \ge 0$ and asymptotically
converges to $\Sigma(\v \omega_0)$. Furthermore,  the set $E(\v \omega_0)$ of relative equilibria with rotation vector $\v \omega_0$
is asymptotically stable in shape space  and every other positive limit set is unstable. \qed  
\end{theorem}

In steady state, the particle motion is characterized by a constant (consensus) twist %This set is characterized by a constant twist
 $\v \xi_0=[\v v_0^T,\v \omega_0^T]^T$.
%\[ \v \xi_0= \left[\begin{array}{c}
%\v s_0\\
%\v \omega
%\end{array}\right].
%\]
The corresponding screw parameters are a pitch $p_0 = <\v v_0,\v \omega_0>/\vn{\v \omega_0}^2$, an axis $l_0 = \{\v v_0 \times \v
\omega_0/\vn{\v \omega_0}^2 + \lambda \v \omega_0, \lambda \in \RR\}$ and a magnitude $M_0 = \vn{\v \omega_0}$.
Therefore the control law (\ref{concon2}) stabilizes all the particles to a relative equilibrium whose pitch depends on the initial conditions of the particles. To
reduce the dimension of the equilibrium set we combine the Lyapunov function (\ref{S}) with the potential
\begin{equation}\label{W}
Q(\v x, \v \omega_0) = \frac{N}{2}\left (\frac{<\v \omega_0, \v x_{\av}>}{\vn{\v \omega_0}} - \alpha \right)^2, \quad \quad \alpha \in [0,1),
\end{equation}
that is minimum when all the particles follow a trajectory with the same pitch $p_0 = \alpha$. This leads to the control law
\begin{equation}\label{concon4}
\begin{array}{lcl}
\v u_k &=& R^T_k \Big [ \v \omega_0 + \Big [ \left(\v r_k - \v r_{\av}\right)\times \v \omega_0 - \v x_{\av}  \\ 
&&\hspace{-0.1 cm}\left. \left.+ \Big (\frac{<\v \omega_0, \v x_{\av}>}{\vn{\v \omega_0}}-\alpha\Big)\frac{\v \omega_0}{\vn{\v \omega_0}} \right]\times \v x_k\right],
\end{array}
\end{equation}
for $k=1,\ldots, N$, which guarantees that $Q+S$ is non-increasing along the solutions. 
\begin{theorem}\label{th:2}
Consider model (\ref{model}) with the control law (\ref{concon4}).
The closed-loop vector field is invariant  under an action of
the translation group $\RR^3$ on position variables $r_{k}$. Every solution exists for all $t \ge 0$ and asymptotically
converges to $\Sigma(\v \omega_0)$. Furthermore, the set of relative equilibria with rotation vector $\v \omega_0$
and pitch $\alpha$ is asymptotically stable in shape space
and every other positive limit set is unstable. \qed
\end{theorem}
%\begin{theorem}\label{th:2}
%Consider model (\ref{model}) equipped with the control law (\ref{concon4}).
%All the solutions converge to the set of equilibria of 
%\begin{equation}\label{s02}
% \dot{\v s}_k =-\Pi_{{\v x_k}}\left(
%\v s_k -\v v^a_{\av} + \left(\frac{<\v \omega, \v x_{\av}>}{\vn{\v \omega}}-\alpha\right)\frac{\v \omega}{\vn{\v \omega}}\right).
%\end{equation}
%The only asymptotically stable equilibria of (\ref{s02}) are relative equilibria
%of (\ref{model}) characterized by an axis of rotation (parallel to $\v \omega$), a pitch $p_0 =\alpha$ and a magnitude $\vn{\v \omega}$. Every other equilibrium of (\ref{s02}) is unstable. \qed
%\end{theorem}
%The proof is reported in Appendix C.

The control law (\ref{concon4}) stabilizes all the particles to a relative equilibrium whose magnitude and pitch are fixed by
the design parameters $\alpha$ and $\vn{\v \omega_0}$. In particular, acting on $\alpha$ it is possible to separate circular
 relative equilibria ($\alpha = 0$) from helical relative equilibria ($\alpha \in (0,1)$). 

It is worth noting that when $\v \omega_0$ is set to zero the control law (\ref{concon2}) reduces to
\begin{equation}\label{conconpar}
\v u_k = R^T_k \left(\v x_k \times \v x_{\av} \right),  \quad \quad k=1,\ldots, N.
\end{equation}
This control law stabilizes parallel formations and has been studied in Section \ref{sec:parallel}.

%The control law (\ref{concon4}) stabilizes the relative  of one particle relying on the vectors $\v \omega$.
\subsection{Dynamic shape control laws for stabilization of screw formations}\label{sec:consome}
Because the control laws (\ref{concon2}) and (\ref{concon4}) depend on the vector $\v \omega_0$, the resulting closed-loop vector
field is not invariant under an action of the rotation group $SO(3)$ on the rotation vaiables. An important consequence is that additional information is
required besides the \emph{relative} configurations among the particles. To overcome this obstacle we propose a consensus approach to
reach an agreement about the direction of the axis of rotation. We provide each particle with a consensus variable $\v \omega_k$, and we denote by $\v \omega^a_k = R_k \v \omega_k$ the same quantity expressed in a (common) spatial reference frame.
The potential
\begin{equation}\label{L_omega}
U(\v \omega^a) = \frac{N}{2}\sum_{k=1}^N \vn{\v \omega^a_k-\v \omega^a_{\av}}^2,
\end{equation}
where $\v \omega^a$ is the stacking vector of the vectors $\v \omega^a_1\,\ldots,\omega^a_N$, decreases along the gradient
dynamics
\begin{equation}\label{consenso}
\dot{\v \omega}_k^a = \sum_{j=1}^N \left(\v \omega^a_j - \v \omega^a_k\right),  \quad \quad k=1,\ldots, N.
\end{equation}
Expressing (\ref{consenso}) in the body reference frame we obtain
\begin{equation}\label{consenso2}
\dot{\v \omega}_k = \hat{\v u}^T_k \v \omega_k+ \sum_{j=1}^N R_k^T R_j \v \omega_j - \v \omega_k,
\end{equation}
for $k=1,\ldots, N,$ and we observe that the dynamics (\ref{consenso2}) are invariant
under an action of the symmetry group $SO(3)$. It turns out that
the dynamic control law resulting from the coupling between the
consensus dynamics (\ref{consenso2}) with the control law
(\ref{concon2}) leads to the shape
control law
\begin{equation}\label{concon5}
\begin{array}{lll}
\hspace{-2 mm}\v u_k &=&  \v \omega_k + \left[R^T_k\left(\v r_k - \v r_{\av}\right)\times \v \omega_k -R^T_k \v x_{\av}\right] \times \v e_1,\\
\hspace{-2 mm}\dot{\v \omega}_k &=& \hat{\v u}^T_k \v \omega_k+ \sum_{j=1}^N R_k^T R_j \v \omega_j - \v \omega_k,
\end{array}
\end{equation}
for $k=1,\ldots, N$.
In the sequel, we denote by 
\[C_{\v \omega} =\{\v \omega_{k}^{a}, k=1,2,\ldots, N\,: \v \omega^{a}_{k} =\v \omega^{a}_{j}, j=1,2,\ldots,N\}
\]
the set of consensus states for the controller variables $\v \omega^{a}_{k}, k=1,2,\ldots,N$~\footnote{From here on we will denote with $C_{\v \eta}$ the set of consensus states for the variables $\v \eta_k^{a}\in \RR^{3}, k=1,2,\ldots,N$.}.  
\begin{theorem}\label{th:final}
Consider model (\ref{model}) with the dynamic control law (\ref{concon5})
%\begin{equation}\label{concon5}
%\begin{array}{rcl}
%\v u_k &=&  \omega_0 \frac{\v \omega_k}{\vn{\v \omega_k}} + \omega_0\left[R^T_k\left(\v r_k - \bar{\v
%r}\right)\times \frac{\v \omega_k}{\vn{\v \omega_k}} -R^T_k\bar{\v x}  \right.\\
%&&\left.+ \omega_0(<\frac{\v \omega_k}{\vn{\v \omega_k}}, \v
%e_1>-\alpha) \frac{\v \omega_k}{\vn{\v \omega_k}}\right] \times \v e_1,\\
%\dot{\v \omega}_k &=& \hat{\v u}^T_k \v \omega_k+ \sum_{j=1}^N R_k^T R_j \v \omega_j - \v \omega_k.
%\end{array}
%\end{equation}
%where $\omega_0 \in \RR$,$\,\alpha \in [0,1)$ and
The closed-loop vector field is invariant under an action of the group $SE(3)$ on the state variables $(\v r_k,R_k)$ and an action of the group $\RR^3$ on the controller variables $\v \omega_k^{a}$. Every solution exists for all $\,t \ge 0$, and asymptotically converges to $\Sigma \times C_{\v \omega}$. Furthermore, $E \times C_{\v \omega}$ is asymptotically stable in the (extended) shape space and every other positive limit set is unstable.\qed
\end{theorem}
%The Theorem is proved in Appendix D.

\begin{figure}[t]
      \centering
      \begin{tabular}{c c}
      \hspace{-1.2 cm}
    \includegraphics[scale=0.30]{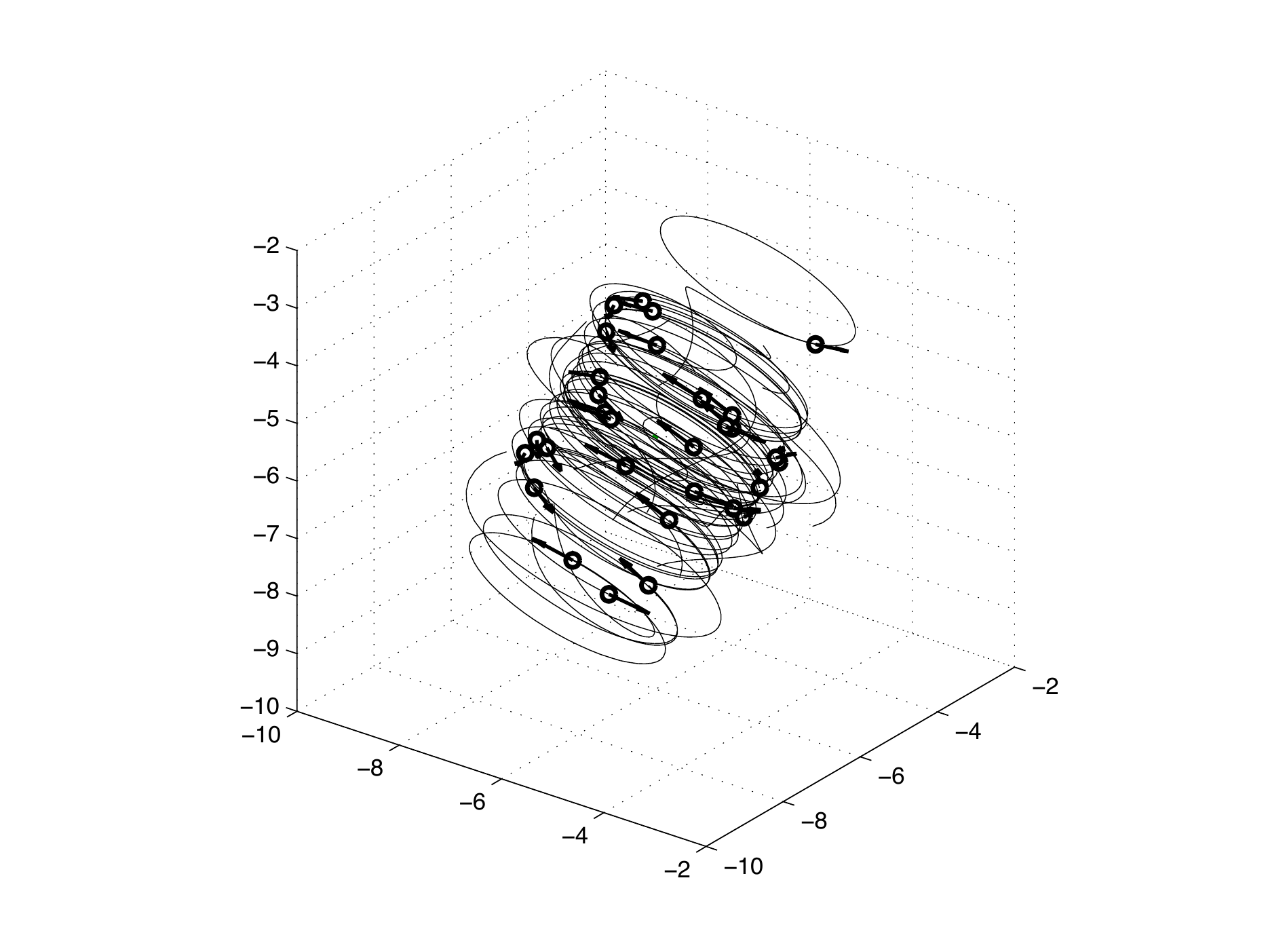}&\hspace{-1.8 cm} \includegraphics[scale=0.32]{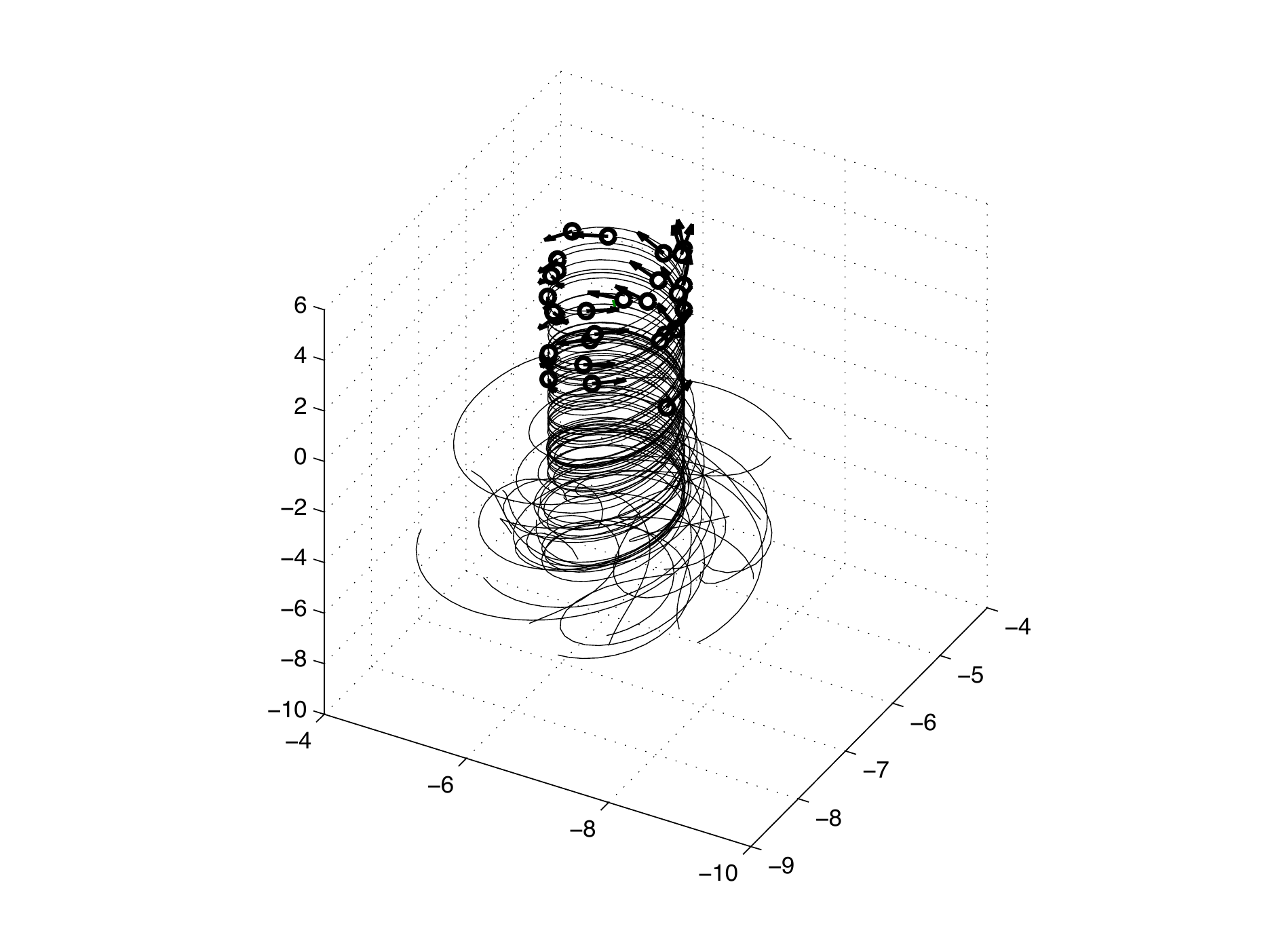}
      \end{tabular}
      \caption{Relative equilibria stabilized with control law (\ref{concon5})}
      \label{fig:sim2}
   \end{figure}

\begin{remark}\label{rem2}
The control law (\ref{concon5}) is the ``dynamic'' version of the control law (\ref{concon2}) and therefore stabilizes all the
particles to a relative equilibrium with arbitrary pitch.
%In Theorem \ref{th:final} we
%recover the convergence properties of the control law (\ref{concon3}) analyzed in Theorem \ref{th:1} without relying on a fixed
%vector $\v \omega$.
To assign to the pitch a desired value it is sufficient to derive the dynamic version of (\ref{concon4}) where consensus dynamics determine a common $\v \omega_0$.
%It is worth noting that the same can be with respect to the control law (\ref{concon4}).
%In this case the only asymptotically stable set is the set of relative equilibria of (\ref{model}) where all the particles move
%about the same axis with desired
%pitch and magnitude.
\end{remark}
In Fig. \ref{fig:sim2} are depicted circular and helical formations stabilized by means of the control law  (\ref{concon5}).

\subsection{Stabilization to a specific screw motion: symmetry breaking} \label{sec:break}
In several applications like sensor networks or formation control, it can be of particular interest to stabilize the motion to a desired screw. To do so, we must break the symmetry
of the control laws presented in the preceding sections. From Section \ref{sec:all-to-all1} we know that a screw is encoded by a constant six dimensional vector $\v \xi_{0} = [\v v_{0}^{T},\v \omega_{0}^{T}]^{T}$.

Consider a virtual particle with dynamics
\begin{equation}\label{model_ref}
\begin{array}{r c l}
\dot{\v r}_0&=&\v x_0\\
\dot{\v x}_{0} &=& \v \omega_{0} \times \v x_{0}. 
\end{array}
\end{equation}  
The particle describes a screw motion characterized by a magnitude  $M_{0} = \vn{\v \omega_{0}}$, an axis $l_{0} = \frac{1}{M_{0}^{2}}\v \omega_{0} \times \left(\v v_{0} \times \v \omega_{0}\right)   + \lambda \v \omega_{0}$ and a pitch $p_{0} = \frac{1}{M_{0}^{2}}<\v x_{0}, \v \omega_{0}>$, where $\lambda \in \RR$ and $\v v_{0} = \v x_{0} + \v r_{0} \times \v \omega_{0}$.   
In the case in which all the particles receive information from the virtual particle, the control law (\ref{concon2}) can be modified as   
\begin{equation}\label{u_ref}
\v u_k = R^T_k \left[\v \omega_{0} + \left( \v v_k^a-\tilde{\v v}^a_{\av}\right)\times \v x_{k}\right].
\end{equation}
for $k=1,2,\ldots,N$,  where $\tilde{\v v}^a_{\av} =\frac{1}{N+1} \sum_{j = 0}^{N}\v v_j^a$. \salt

\begin{proposition}\label{prop:break}
Consider the closed-loop system given by (\ref{model}) and the control law (\ref{u_ref}).
Every solution exists for all $t \ge 0$ and asymptotically
converges to $\Sigma(\v \omega_{0})$. Furthermore,  the set of relative equilibria with rotation vector $\v \omega_{0}$, pitch $p_{0}$ and axis $ l_{0}$ is asymptotically stable and every other positive limit set is unstable.
 \qed
\end{proposition}
%The proof is reported in Appendix E.
\begin{figure}[t]
      \centering
  %    \begin{tabular}{c c}
%\includegraphics[scale=0.25]{helix_break.eps} & 
 \includegraphics[scale=0.9]{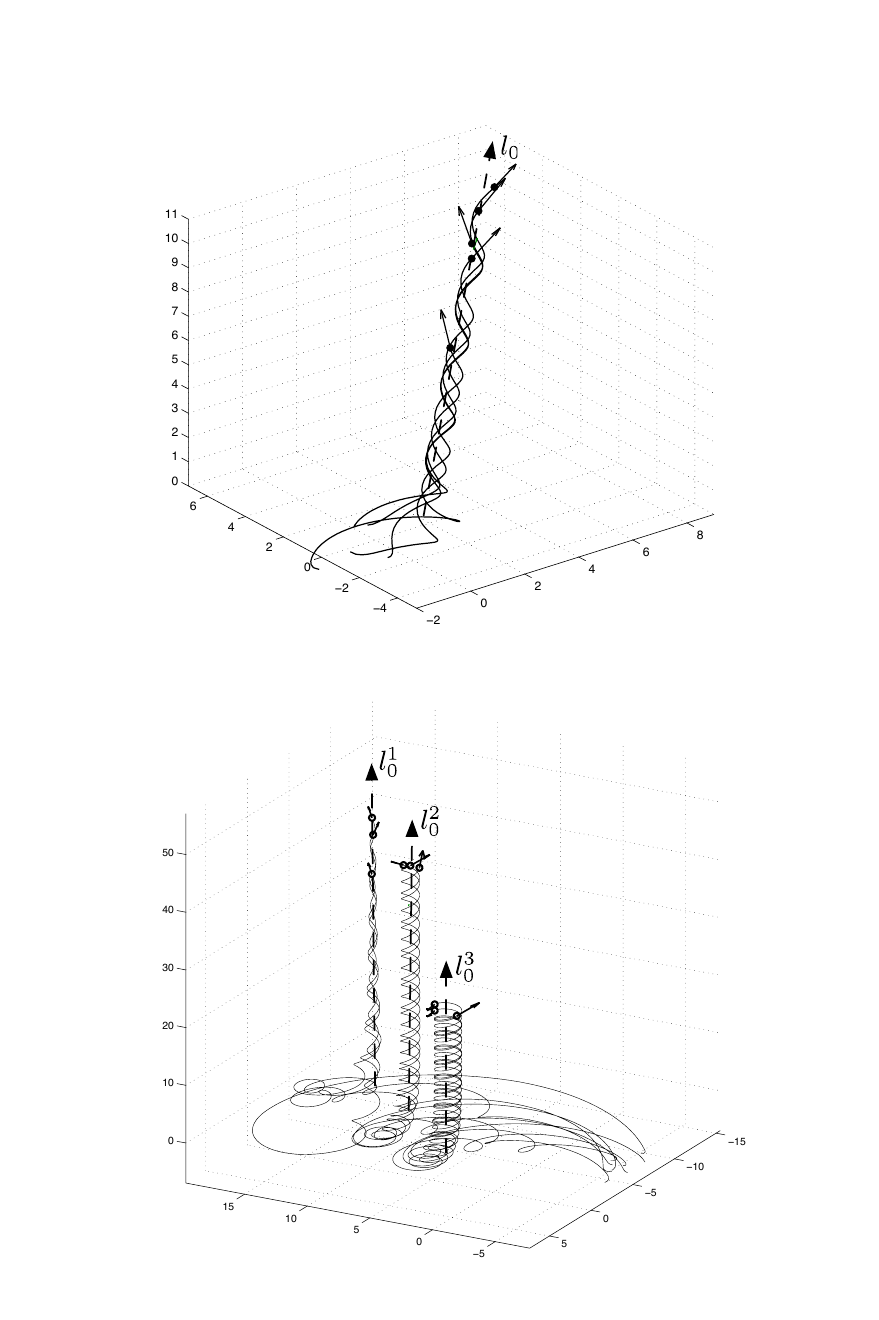}
%     %\hspace{-0.5 cm}  (a)\\
%     \hspace{-0.5 cm} \includegraphics[scale=0.50]{helix_lim2.eps}\\
%     % \hspace{-0.5 cm}  (b)
% %\hspace{-1.8 cm}\includegraphics[scale=0.4]{circ.eps}&\hspace{-2.8 cm} \includegraphics[scale=0.4]{helix.eps}
   %   \end{tabular}
      \caption{On the top: Helical formation stabilized with the control law (\ref{u_ref}). The parameters of the helix are set to $p_{0} = 0.5$, $\v \omega_{0} = [1,1,1]^{T}$ and $l_{0} = [1,-1,0]^{T} + \lambda \v \omega_{0}$. On the bottom: Helical formations stabilized with the control law (\ref{u_ref2}). Each subgroup converges to a different screw defined by a different axis and a different pitch.}  %In (a) $\alpha=0$, while in (b) $\alpha=0.5$. The parameter $\omega_0$ is set to
%      %$1$.
%      The communication graph ${\mathcal G}(t)$  is time varying and uniformly connected.
%      }
      \label{fig:break}
   \end{figure}
This approach is well suited to stabilize subgroups of particles to different screw formations. To this end is sufficient to define a virtual particle for each subgroup and to fix the parameters of the desired screw motions. 
Consider $M$ subgroups of particles $B_{1}\ldots, B_{M}$. For simplicity let the cardinality of each group  be $n$. Define $n$ virtual particles obeying the following dynamics:
\begin{equation}\label{model_ref2}
\begin{array}{r c l}
\dot{\v r}^{i}_0&=&\v x_0^{i}\\
\dot{\v x}^{i}_{0} &=& \v \omega_{0}^{i} \times \v x_{0}^{i}. 
\end{array}
\end{equation}  
 for $i=1,\ldots,M$. Define  $\tilde{\v v}^{i}_{\av} =\frac{1}{n+1} \left(\sum_{j \in B_{i}}\v v_j^a + \v v^{i}_{0}\right)$, where $\v v^{i}_{0} = \v x_{0}^{i} + \v r_{0}^{i} \times \v \omega_{0}^{i}$ and $\v v_j^a = \v x_j + \v r_{j} \times \v \omega_{0}^{i},\, j\in B_{i}$  (where, with a little abuse of notation, we dropped the apex $a$ in the average velocity).     

The following control law generalizes (\ref{u_ref}): 
\begin{equation}\label{u_ref2}
\v u_k = R^T_k \left[\v \omega^{i}_{0} + \left( \v v_k^a-\tilde{\v v}^{i}_{\av}\right)\times \v x_{k}\right],\quad \quad k \in B_{i}
\end{equation}
for $i=1,\ldots,M$.

As a direct corollary of Proposition \ref{prop:break} the control law (\ref{u_ref2}) stabilizes the particles in each group $B_{i}, i=1,\ldots,M$ to a screw motion defined by $\v \xi_{0}^{i} = [\v {v_{0}^{i}}^{T}, {\v \omega_{0}^{i}}^{T} ]^{T}$.
In Fig. \ref{fig:break} different motion patterns, obtained by adopting control laws (\ref{u_ref}) and (\ref{u_ref2}), are displayed.    

All the control laws presented until this point stabilize the relative
equilibria of (\ref{model}) under the assumption of all-to-all
communication among the particles. In Section \ref{sec:shape} we relax
this requirement by substituting the quantities in (\ref{concon5})
that require global information with consensus variables obeying consensus dynamics.

Before detailing the approach, in the following section we review
some concepts about consensus in Euclidean space
and we summarize some graph theoretic notions that are needed to address the problem
in a limited communication setting.

\section{Communication graphs and consensus
dynamics in Euclidean space}\label{sec:consensus} In this section we review some recent results on the consensus problem.
Consider a group of agents with limited communication capabilities; in this context it is useful to describe the communication
topology by using the notion of \emph{communication graph}.

Let $G=({\mc V},{\mc E},A)$ be a weighted digraph (directed graph)
where ${\mc V}=\{v_1,\ldots,v_N\}$ is the set of nodes, ${\mc
E}\subseteq {\mc V} \times {\mc V}$ is the set of edges, and ${A}$
is a weighted adjacency matrix with nonnegative elements $a_{kj}$.
We assume that there are no self-cycles
i.e. $a_{kk}=0,\, k=1,2,\ldots, N$.

The graph Laplacian $L$ associated to the graph $G$ is defined as
\[
L_{kj}=\left\{
\begin{array}{ll}
\sum_i a_{ki}, &   j=k\\
-a_{kj},       &   j\ne k.
\end{array}
\right.
\]
%The digraph $G$ is called {\em strongly connected} if
%and only if any two distinct nodes of the graph can be connected via a path that follows the direction of the edges of the
%digraph.
The $k$-th row of $L$ is defined by $L_k$. The in-degree
(respectively out-degree) of node $v_k$ is defined as
$d_k^{in}=\sum_{j=1}^N a_{kj}$ (respectively
$d_k^{out}=\sum_{j=1}^N a_{jk}$). The digraph $G$ is said to be
{\em balanced} if the in-degree and the out-degree of each node
are equal, that is,
\[
\sum_j a_{kj} = \sum_j a_{jk}, \quad \quad k=1,2,\ldots, N.
\]
%It is both of theoretical and practical interest to consider
%time-varying communication topologies. For example, in a network
%of moving agents, some of the existing links can fail and new
%links can appear when other agents enter an effective range of
%detection.
If the communication topology is time varying, it can be described
by the time-varying graph $G(t)=({\mc V},{\mc E}(t),A(t))$, where
$A(t)$ is piece wise continuous and bounded and $a_{kj}(t)\in
\{0\}\cup[\eta,\gamma], \forall\, k, j,$ for some finite scalars
$0<\eta\leq\gamma$ and for all $t\geq 0$. The set of neighbors of
node $v_k$ at time $t$ is denoted by ${\mc N}_k(t) \triangleq
\{v_j \in {\mc V} : a_{kj}(t)\geq\eta\}$. We recall two
definitions that characterize the concept of uniform connectivity
for time-varying graphs.
\begin{definition}
Consider a graph $G(t)=({\mc V},{\mc E}(t),A(t))$. A node $v_k$ is
said to be connected to node $v_j$ ($v_j \ne v_i$) in the interval
$I=[t_a,t_b]$ if there is a path from $v_k$ to $v_j$ which
respects the orientation of the edges for the directed graph
$({\mc N},\cup_{t \in I} {\mc E}(t),\int_I A(\tau)d \tau)$.
\end{definition}
\begin{definition}
$G(t)$ is said to be uniformly connected if there exists a time
horizon $T>0$ and an index $k$ such that for all $t$ all the nodes
$v_j$ ($j\neq k$) are connected to node $k$ across $[t, t+T]$.
\end{definition}
Consider a group of $N$ agents with state $\v p_k \in P$, where
$P$ is an Euclidean space. The communication between the
$N$-agents is defined by the graph $G$: each agent can sense only
the neighboring agents, i.e. agent $j$ receives information
from agent $i$ if and only if $i\in {\mc N}_j(t)$. %We use the notation $k \sim j$ to indicate the presence of a communication link from
%agent $j$ to agent $k$, i.e. $\,k \sim j\,$ iff $v_j\in {\mc N}_k$.

Consider the continuous dynamics \ee \label{cons1} \dot{\v
p}_k=\sum_{j=1}^{N} a_{kj}(t)(\v p_j-\v p_k), \quad \quad 
k=1,2,\ldots, N.\eee Using the Laplacian definition, (\ref{cons1})
can be equivalently expressed as \ee \label{cons2} \dot{\v
p}=-\tilde{L}(t)\,\v p, \eee where $\tilde L = L \otimes I_3$ and
$\v p = (\v p_1^T, \ldots, \v p_N^T)^T$. Algorithm (\ref{cons2})
has been widely studied in the literature and asymptotic
convergence to a consensus value holds under mild assumptions on
the communication topology.
%In \cite{OlMu} some conditions on the communication topology have
%been given, under which the solutions of \neq{cons2} converge to $\bar{x}(0) 1$, where $\bar{x}(0)$ is the average of the
%initial condition, i.e. $\bar{x}(0)=\sum_{k \in {\mc I}} x_k(0)$. The authors showed that, if the communication graph switches
%among a finite number of strongly connected and balanced graphs then the consensus on the average value is asymptotically
%reached. The strongly connection hypothesis can be relaxed in light of the results of \cite{Mo}.
The following theorem summarizes some of the main results in
\cite{Mocdc}, \cite{Mo} and \cite{OlMu}.
\begin{theorem}\label{ConsEucl}
Let $P$ be a finite-dimensional Euclidean space. Let $G(t)$ be a
uniformly connected digraph and $L(t)$ the corresponding
Laplacian matrix bounded and piecewise continuous in time. %Consider the linear systems \neq{cons2} and \neq{cons2d}
%\ee \label{cons3} \dot{x}=-L(t)x,\quad \quad x\in X^N, \eee
The solutions of (\ref{cons2}) asymptotically converge to a consensus value $\beta {\bf 1}$ for some $\beta \in P$. Furthermore
if $G(t)$ is balanced for all $t$, then $\beta =\frac{1}{N} \sum_{i=1}^{N}\v p_i(0)$.\qed
\end{theorem}
A general proof for Theorem \ref{ConsEucl} is based on the
property that the convex hull of vectors $\v p_k \in W$ is non
expanding along the solutions. For this reason, the assumption
that $W$ is an Euclidean space is essential (see e.g. \cite{Mo}).
Under the additional balancing assumption on $G(t)$, it follows
that ${\bf 1}^T L(t) = 0$, which implies that the average
$\frac{1}{N} \sum_{j \in {\mc I}} \v p_j $ is an invariant
quantity along the solutions.

\section{Stabilization of relative equilibria in the presence of limited communication}\label{sec:shape}
Consider the control laws (\ref{cont1}) and (\ref{concon5}). By
following the approach presented in \cite{scardovi:2007}
%\begin{equation*}
%\begin{array}{rcl}
%\v u_k &=&  \v \omega_k + \left[R^T_k\left(\v r_k - \bar{\v
%r}\right)\times \v \omega_k -R^T_k\bar{\v x}  \right.\\
%&&\left.+ (<\v \omega_k, \v
%e_1>-\alpha)\v \omega_k\right] \times \v e_1,\\
%\dot{\v \omega}_k &=& R_k^T \hat{\v u}_k + \sum_{j=1}^N R_k^T R_j
%\v \omega_j - \v \omega_k.
%\end{array}
%\end{equation*}
we substitute the quantities that require all-to-all
communication, i.e. $\v r_{\av}$ and $\v x_{\av}$, by local
consensus variables. This leads to a generalization of the control
laws (\ref{cont1}) and (\ref{concon5}) to uniformly connected
communication graphs. We consider first the problem of stabilizing a
parallel formation.
\subsection{Stabilization of parallel formations with limited communication}
We replace the control law (\ref{cont1})
%\begin{equation}
%\begin{array}{rcl}
%q_k&=&<{\v x}_{\av}, \v y_k> \\
%h_k&=&<{\v x}_{\av}, \v z_k>, \quad \quad k=1,\ldots,N,
%\end{array}
%\end{equation}
with the local control law
\begin{equation}\label{cont_limited}
\v u_k=R_k^T(\v x_k \times \v b_k^{a}), \quad \quad k=1,\ldots, N,
%\begin{array}{rcl}
%q_k&=&<\v b_k, \v y_k> \\
%h_k&=&<\v b_k, \v z_k>, \quad \quad k=1,\ldots,N,
%\end{array}
\end{equation}
where $\v b_k^a$ is a consensus variable obeying the consensus
dynamics
\begin{equation}\label{cons_w}
\dot{\v b}^{a}_k = -\sum_{j=1}^N L_{kj} \v b^{a}_j, \quad \quad k=1,\ldots, N,
\end{equation}
with arbitrary initial conditions $\v b^{a}_k(0), k=1,\ldots,N$.
Before detailing the convergence analysis we express
(\ref{cont_limited}) and (\ref{cons_w}) in shape coordinates by moving to a local reference frame. 
Then (\ref{cont_limited}) rewrites as
\begin{equation}\label{cont_limited_shape}
%\begin{array}{rcl}
%q_k&=&<\v b_k, \v e_2> \\
%h_k&=&<\v b_k, \v e_3>, \quad \quad k=1,\ldots,N,
%\end{array}
\v u_k=(\v e_1 \times \v b_k),\quad \quad  k=1,\ldots, N,
\end{equation}
and (\ref{cons_w}) as
\begin{equation}\label{cons_w_shape}
\dot{\v b}_k = \hat{u}_k^T \v b_k - \sum_{j=1}^N L_{kj} R_k^T R_j
\v b_j,
\end{equation}
where $\v b_k(0) = R^T_k(0) \v b^{a}_k(0)$, $k=1,\ldots, N$. 
The following result characterizes the convergence properties of the resulting closed-loop system.
\begin{theorem}
Consider model (\ref{model}) with the control law (\ref{cont_limited_shape}),(\ref{cons_w_shape}).  The closed-loop vector field is invariant  under an action of the group $SE(3)$ on the state variables
$(\v r_k,R_k)$ and  an action of the group $\RR^3$ on the consensus variables $\v b^{a}_k$. Suppose that the communication graph $G(t)$ is uniformly connected and that $L(t)$ is bounded and piecewise continuous.
Then every solution exists for all $t \ge 0$
and asymptotically converge to $\Sigma(0)\times C_{\v b}$. Furthermore, the set
$E(0)\times C_{\v b}$ is asymptotically stable in the (extended) shape space and every other positive limit set is unstable. \qed
\end{theorem}
%Consider now the problem of balancing the velocity vectors. We
%replace
%\begin{equation}
%\begin{array}{rcl}
%q_k&=&-<{\v x}_{\av}, \v y_k> \\
%h_k&=&-<{\v x}_{\av}, \v z_k>, \quad \quad k=1,\ldots,N,
%\end{array}
%\end{equation}
%with
%\begin{equation}\label{cont_limited2}
%\begin{array}{rcl}
%q_k&=&-<\v p_k, \v y_k> \\
%h_k&=&-<\v p_k, \v z_k>, \quad \quad k=1,\ldots,N,
%\end{array}
%\end{equation}
%where $\v p_k$ is a consensus variable obeying the consensus
%dynamics
%\begin{equation}
%\dot{\v p}_k = -\sum_{j=1}^N L_{kj} \v p_j + \frac{d}{dt} \v x_k,
%\end{equation}
%where $\v p_k(0) = \v x_k(0)$ for every $k$.

%In shape coordinates we obtain

%\begin{equation}\label{cont_limited_shape2}
%\begin{array}{rcl}
%q_k&=&-<\v b_k, \v e_2> \\
%h_k&=&-<\v b_k, \v e_3>, \quad \quad k=1,\ldots,N,
%\end{array}
%\end{equation}
%where $\v b_k$ is a consensus variable obeying the consensus
%dynamics
%\begin{equation}
%\dot{\v b}_k = \hat u_k^T\left(\v b_k - \v e_1\right)-\sum_{j=1}^N
%L_{kj} R_k^T R_j \v p_j,
%\end{equation}
%where $\v b_k(0) = \v e_1$ for every $k$.\salt
%\begin{theorem}
%Suppose that the communication graph $G(t)$ is uniformly connected
%and balanced for all $t\geq0$ and that $L(t)$ is bounded and
%piecewise continuous. Then all the solutions of the system
%(\ref{model3}) with the control (\ref{cont_limited2})
%asymptotically converge to an equilibrium. Moreover, the only
%stable limit set is the set of balanced states characterized by
%${\v x}_{\av}=0$.\qed
%\end{theorem}
%The proof is reported in Appendix G.

\subsection{Stabilization of screw formations in the presence of limited communication}
We finally address the problem of stabilizing screw relative
equilibria in the presence of limited communication. The procedure to generalize the control law
(\ref{concon5}) is the same as
outlined in the previous section and therefore is omitted. Consider the dynamic control law 
\begin{equation}\label{concon6}
\begin{array}{rcl}
\v u_k &=& \hspace{-0.1 cm} \v \omega_k + \left(\v \omega_k \times \v c_k -\v b_k \right) \times \v e_1\\
\dot{\v \omega}_k &=&\hspace{-0.1 cm} \hat{\v u}_k^T \v \omega_k -
\sum_{j=1}^N
L_{kj} R_k^T R_j \v \omega_j\\
\dot{\v b}_k &=&\hspace{-0.1 cm}\hat{\v u}_k^T \v b_k -
\sum_{j=1}^N L_{kj}R_k^T R_j
\v b_j\\
\dot{\v c}_k &=& \hspace{-0.1 cm}\hat{\v u}_k^T \v c_k - \v e_1
\hspace{-0.1 cm} -\sum_{j=1}^N  L_{kj}R_k^T R_j \v c_j -
\hspace{-0.1 cm}\sum_{j=1}^N L_{kj}R_k^T \v r_j,
\end{array}
\end{equation}
for $k=1,\ldots, N$, and define $\v \omega_k^a=R_k \v \omega_k,\, \v b_k^a=R_k \v b_k,\, \v
c_k^a=R_k \v c_k + \v r_k$.
\begin{theorem}\label{th:final2}
 Consider model (\ref{model}) with the control law (\ref{concon6}).
The closed-loop vector field is invariant  under an action of the group $SE(3)$ on the state variables
$(\v r_k,R_k)$ and  an action of the group $\RR^3 \times \RR^3 \times \RR^3$ on the consensus variables $(\v \omega^{a}_k,\v b^{a}_k,\v c^{a}_k)$.  Suppose that the communication graph $G(t)$ is uniformly connected and that $L(t)$ is bounded and piecewise continuous.
%The consensus variables
%$\v \omega^{a}_k$, $\v b^{a}_k$, and $\v c^{a}_k$ exponentially converge to consensus values $\v \omega_{0}$, $\v b_{0}$, and
%$\v c_{0}$, respectively. 
Then every solution exists for all $t \ge 0$
and asymptotically converge to $\Sigma\times C_{\v \omega}\times C_{\v b}\times C_{\v c}$. Furthermore, the set
$E\times C_{\v \omega}\times C_{\v b}\times C_{\v c}$ is asymptotically stable in the (extended) shape space and every other positive limit set is unstable. 
\end{theorem}

%\begin{theorem}\label{th:final2}
%Let $G(t)$ be a directed and uniformly connected communication
%graph and $L(t)$ the corresponding bounded and piecewise
%continuous Laplacian matrix. Consider model (\ref{model}) equipped
%with the control law
%where the initial conditions of the consensus variables are randomly selected. The resulting closed-loop vector field is
%invariant under an action of $SE(3)$. The only asymptotically stable set in the shape space is the set of relative equilibria of
%model (\ref{model}). Every other equilibrium is unstable.\qed
%\end{theorem} 
%The proof is reported in Appendix H. 

It is important to note that the control law (\ref{concon6}) does
not require all-to-all communication among the particles. In
particular the convergence properties of Theorem \ref{th:final}
are here recovered in the presence of limited communication, for
directed, time-varying (but uniformly connected) communication
topologies. Furthermore, following the approach proposed in \cite{SePaLe_limited}, it is possible to extend the symmetry-breaking approach presented in Section \ref{sec:break} to the limited communication scenario. This can be done redefining the graph Laplacian by adding a directed link connecting every particle to a virtual particle. The uniformly connectedness assumption on the new graph guarantees convergence to the desired screw motion.   

Due to space constraints we do not report here the details, the interested reader is refereed
to \cite{SePaLe_limited} where the planar case is considered.     

\section{Discussion on possible applications}\label{sec:applications}       
In this paper models of point-mass particles at constant speed are considered.
From the engineering and application-oriented  prospective, they are a strong simplification of 
the dynamic models that can be used in ``real world'' applications. To introduce more sophisticated models in our scheme, a reasonable solution is to decouple the collective design problem (that we have addressed in the present paper) with a trajectory tracking problem where the details about the system dynamics are taken in account. This means that each vehicle is provided with a trajectory ``planner'' that designs the required trajectory by exchanging information with the other vehicles. A second module, namely a tracking controller, must be designed to ensure that the discrepancy between the actual trajectory and the designed one is minimized. This module incorporates the details about the dynamics of the system and is completely decoupled from the other vehicles.

A particularly interesting application is the collection of sensor data with underwater gliders. Underwater gliders are autonomous vehicles that
rely on changes in vehicle buoyancy and internal mass redistribution for regulating
their motion. They do not carry thrusters or propellers and have limited external
moving control surfaces. For these vehicles only a subset of the relative equilibria
may be realized, and they correspond to motion (at constant speed) along circular
helices and straight lines \cite{phdbhatta}. In particular, for equilibrium
motion along a circular helix, the axis of the helix must be aligned with the direction
of gravity. This suggests to apply the control laws presented in the present paper, fixing the direction of the rotation axis to $\v \omega_0 = -c \,[0, 0, 1]^{T}$, where $c$ is a constant positive scalar, to plan the desired trajectories. The parameters of the desired helical motion, and consequently of the control law of the planner,  can be chosen on the basis of energy efficiency criteria (which depend on the glider's parameters) and to concentrate the data collection at the desired location. The problem of designing a trajectory tracking controller for underwater gliders has been addressed in \cite{phdbhatta} and is beyond the scope of the present work.             
\section{Conclusions}
We propose a methodology to stabilize relative equilibria in a
model of identical, steered particles moving in three-dimensional
Euclidean space. Observing that the relative equilibria can be
characterized by suitable invariant quantities, we formulate the
stabilization problem as a consensus problem. The formulation leads to
a natural choice for the Lyapunov functions. Dynamic control laws
are derived to stabilize relative equilibria in the presence of
all-to-all communication and are generalized to deal with
unidirectional and time-varying communication topologies. It is of
interest (in particular from the application point of view) to
study in the future how to reduce the dimension of the equilibrium
set by breaking the symmetry of the proposed control laws.
%\vspace{3 cm}

\appendix{\bf Appendix}
\section{Proof of Theorem 1}
Since the control law (\ref{cont1}) is independent from the relative spacing of the particles, we can limit our analysis to the reduced dynamics (\ref{model3}).
Plugging (\ref{cont1}) into (\ref{V_x_dot}) yields
\[
\dot{V}= -\sum_{k=1}^N \vn{\v x_k \times \v x_{\av}}^2 \leq 0.
\]
$V$ is positive definite (in the reduced shape space) and non increasing. By the La Salle invariance principle, the solutions of (\ref{model3}) converge to
the largest invariance set where
\begin{equation}\label{set1}
\v x_k \times \v x_{\av}=0, \quad \quad k=1,\ldots,N.
%\begin{array}{rcl}
%<{\v x}_{\av}, \v y_k>&=&0 \\
%<{\v x}_{\av}, \v z_k>&=&0, \quad \quad k=1,\ldots,N.
%\end{array}
\end{equation}
This set is contained in $\Sigma(0)$.
The points where ${\v x}_{\av} = 0$, are global maxima of $V$. As
a consequence this set is unstable. From (\ref{set1}), equilibria where ${\v x}_{\av} \neq
0$ are characterized by the vectors $\v x_k,\, k=1,\ldots,N,$ all
parallel to the constant vector with ${\v x}_{\av}$. Note that this
configuration involves $N-M$ velocity vectors aligned to ${\v
x}_{\av}$ and $M$ velocity vectors anti-aligned with ${\v x}_{\av}$,
where $0\leq M<\frac{N}{2}$. At those points, $\vn{{\v
x}_{\av}}=1-\frac{2M}{N}>\frac{1}{N}$. When $M=0$ we recover the
set of synchronized states (global minima of $V$) which is stable. Every other value of $M$
corresponds to a saddle point (isolated in the shape space) and is therefore unstable. To see this we express $\v x_k$ and ${\v
x}_{\av}$ in spherical coordinates,
\[
{\v x}_{\av}= \vn{{\v x}_{\av}}[\cos \Phi \sin \Theta,\sin \Phi \sin \Theta, \cos \Theta]^T, %\quad 0\leq\Theta \leq \pi, \quad
%0\leq\Phi \leq 2\pi,
\]
\[
\v x_k= [\cos \phi_k \sin \theta_k,\sin \phi_k \sin \theta_k, \cos \theta_k]^T, %\quad 0\leq\theta_k \leq \pi, \quad 0\leq\phi_k
%\leq 2\pi.
\]
where $\theta_k,\Theta \in [0,\pi]$ and $\phi_k,\Phi \in
[0,2\pi)$. By expressing $V$ with respect to spherical coordinates we obtain \ee
\label{lypr}
\begin{array}{rcl} 
V&=&\displaystyle \frac{N}{2}\Big(1-\frac{1}{N}\vn{{\v x}_{\av}}\sum_{j=1}^N\big ( \sin \Theta \sin
\theta_j \cos( \Phi- \phi_j) \\
&&+ \cos \Theta \cos \theta_j \big ) \Big). 
\end{array}
\eee
The critical points are characterized by
\[
\v x_k= [\cos \Phi \sin \Theta,\sin \Phi \sin \Theta, \cos
\Theta]^T,\; k=M+1,\ldots,N,
\]
and \[
\begin{array}{rcl}
\v x_k\hspace{-0.1 cm}&=&\hspace{-0.1 cm} [\cos (\Phi + \pi) \sin (\pi-\Theta),\sin (\Phi+\pi) \sin
(\pi-\Theta),\\ 
&&\cos (\pi-\Theta)]^T, \;  k=1,\ldots,M.
\end{array}
 \]

The second derivative of $V$ (with respect to $\theta_j$) is
\[
\frac{\partial^2 V}{\partial \theta^2_j} \hspace{-0.1 cm}=\hspace{-0.1 cm}  \vn{{\v
x}_{\av}} \left(\sin\theta_j \sin\Theta \cos(\Phi-\phi_j)\hspace{-0.7 mm} +\hspace{-0.7 mm} \cos
\Theta \cos \theta_j\right)-\frac{1}{N},
\]
that is positive if $\theta_j= \Theta$ and $\phi_j= \Phi$ and is
negative if $\theta_j= \pi - \Theta$ and $\phi_j= \Phi+\pi$. As a
consequence, a small variation $\delta \theta_j$ at critical
points where $M \neq 0$ increases the value of $V$ if $\theta_j= \Theta$ and
$\phi_j= \Phi$, and decreases the value of $V$ if $\theta_j= \pi
- \Theta$ and $\phi_j= \Phi+\pi$. %When $\Theta\neq\,(0$ mod $\pi$)

We conclude that $E(0)$ (the set of relative equilibria corresponding to parallel motion) is asymptotically stable in the shape space and the other positive limit sets are unstable.
\qedp

\section{Proof of Theorem 2}
$S$ is non negative and, from (\ref{Sdot}), it is non-increasing along the
solutions of (\ref{model}). Then $S$ converges to a limit as $t \rightarrow \infty$.
Furthermore the second derivative $\ddot S$ is bounded (because $\v v_k^{a} -  \v v^a_{\av}$ is bounded for every $k$). From Barbalat's Lemma $\dot S \rightarrow 0$ when $t\rightarrow \infty$ and therefore the solutions converge to the set
$\Gamma$ where
\begin{equation}\label{cond1}
\left(\v v_k^a - \v v^a_{\av}\right)\times \v x_k = 0,
\end{equation}
that characterizes the equilibria of (\ref{s0}). 
Observe that in $\Gamma$, $\dot{\v x}_k = \v \omega_0 \times \v x_k$ and $\v v_k^a$ is constant
for $k=1,\ldots,N$. Therefore $\Gamma \subseteq \Sigma(\v \omega_0)$.  
It remains to prove that the set $E(\v \omega_0)$ is asymptotically stable (in the shape space) and the other sets (in $\Gamma$) are unstable.

We divide the analysis into three parts to analyze $\Gamma$.

i) Suppose that in steady state $\v \omega_0 \times \v x_k \neq 0$
for every $k$. Then (\ref{cond1}) can hold only if $\v v_k^a = \v
v_0$ for every $k$ and for some fixed $\v v_0 \in \RR^3$, this set
defines a global minimum for $S$ and therefore is
asymptotically stable in the shape space. This set corresponds to circular or helical relative equilibria (with axis of rotation parallel to $\v \omega_0$) and is contained in $E(\v \omega_0)$.

ii) Suppose now that in steady state $\dot{\v x}_k = \v \omega_0 \times \v x_k = 0$ for every $k$. From (\ref{cond1}) we obtain
\[
\left(\v v_k^a - \v v^a_{\av}\right)\times \v \omega_0 = 0
\]
for every $k$, which implies $(\v r_k - \v r_{\av}) \times \v
\omega_0 = 0$. Therefore in steady state the Lyapunov function
(\ref{S}) reduces to
\begin{equation}\label{S2}
S= \frac{1}{2}\sum_{k=1}^N \vn{\v x_k - \v x_{\av}}^2.
\end{equation}

This set is characterized by the vectors $\v x_k,\quad k=1,\ldots,N,$ all parallel to the constant vector $\v \omega_0$.
Note that this configuration involves $N-K$ velocity vectors aligned to $\v \omega_0$ and $K$ velocity vectors anti-aligned to $\v
\omega_0$ (or vice-versa), where $0\leq K\leq\frac{N}{2}$. When $K=0$, potential (\ref{S2}) is zero (global minimum), and
therefore the configuration defines an asymptotically stable set. This set corresponds to collinear formations (with
the same direction of motion) parallel to $\v \omega_0$. These configurations are relative equilibria and are contained in $E(\v \omega_0)$.

When $K=\frac{N}{2}$, potential (\ref{S2}) attains a global maximum, and therefore the configuration defines unstable
equilibria. Every other value of $K$ corresponds to a saddle point and is therefore unstable. To see this it is sufficient to
express $\v x_k$ and $\v \omega_0$ in spherical coordinates and to show that
$S$ can decrease under an arbitrary small perturbation (see the proof of Theorem $1$).

iii) It remains to analyze the situation where $\v \omega_0 \times \v x_k \neq 0$ for $k \in G_1$ and $\v \omega_0 \times \v x_j =0$
for $j \in G_2$, where $G_1$ and $G_2$ denote two disjoint groups of particles such that $G_1 \bigcup G_2 = \{1,\ldots,N\}$  and $|G_{1}|=M$ and $|G_{2}|=N-M$. In
such a situation we obtain
\begin{equation}\label{condpro}
\begin{array}{rcl}
\v v_k^a - \v v^a_{\av} &=& 0, \quad \quad k\in G_1\\
(\v v^a_j - \v v^a_{\av})\times \v \omega_0  &=& 0,  \quad \quad j\in G_2,
\end{array}
\end{equation}
where $\v v^a_j\neq \v v^a_{\av},\, j\in G_2$. We call this set $\Lambda$. 
Since
\[
\v v^a_{\av} = \frac{1}{N}\sum_{k\in G_{1}}\v v_k^a +\frac{1}{N}\sum_{j\in G_{2}}\v v_j^a
\]
from (\ref{condpro}) we observe that 
\[
\v v^a_{\av} = \frac{1}{N-M}\sum_{j\in G_{2}}\v v_j^a,
\]
that implies that $(\v r_{k}-\frac{1}{N-M}\sum_{k\in G_{2}}\v r_{k}) \times \v \omega_0=0$ for every $k\in G_{2}$.

Therefore in this set the Lyapunov function (\ref{S}) reduces to
\begin{equation}\label{S3}
\tilde S= \frac{1}{2}\sum_{k\in G_2} \vn{\v x_k-\v x^{G_{2}}_{\av}}^2
\end{equation}
%and attains the value $S=\frac{|G_2|}{2}\vn{\v \omega}^2$. %This set is unstable, to see it take two particles $l,z\in G_2$ and
%perturb the velocity vector $\v x_l$ with a small variation (norm preserving) $\delta$ and the velocity vector $\v x_z$ with
%$-\delta$. Note that $\v v^a_{\av}$ is invariant to this perturbation.
where $\v x^{G_{2}}_{\av} = \frac{1}{N-M}\sum_{k\in G_{2}}\v x_{k}$.
Since $\v v^a_j\neq \v v^a_{\av},\,$ and $\v x_j$ is parallel to $\v
\omega_0$ for every $j\in G_2$, $\v x_{j} \neq \v x_{\av}^{G_{2}}$ for every $j\in G_2$. We conclude from (\ref{S3}) that this set does not correspond to global minima of (\ref{S}). 
Now we prove that this  set is unstable. The first step is to show that this set does not correspond to local minima of (\ref{S}). To this end we express the velocity vectors and the rotation vector in spherical coordinates:
\[
\v x_k= [\cos \phi_k \sin \theta_k,\sin \phi_k \sin \theta_k, \cos \theta_k]^T, %\quad 0\leq\theta_k \leq 
\]
and 
\[
\v \omega_0= [\cos \Phi \sin \Theta,\sin \Phi \sin \Theta, \cos \Theta]^T, %\quad 0\leq\theta_k \leq \pi, 
\]
where $\theta_k,\Theta \in [0,\pi]$ and $\phi_k,\Phi \in
[0,2\pi)$, and we compute the second partial derivative of  (\ref{S}) with respect to a particular direction. Let $\v x_{p}, p\in G_{2}$, be a velocity vector such that $\v x_{p}=-\frac{\v \omega_0}{\vn{\v\omega_0}}$ (notice that such a vector always exists since $\v x_{j} \neq \v x_{\av}^{G_{2}}$ for every $j\in G_2$). We show that the second derivative with respect to $\theta_{p}$ is negative in this set.
After some calculations we arrive at the following expression:
\[
\begin{array}{rcl}
\displaystyle  \frac{\partial^{2}S}{\partial \theta_{p}^{2}} &\displaystyle =& \displaystyle \sum_{k=1}^{N} < \frac{\partial^{2}(\v x_{k}-\v x_{\av})}{\partial \theta_{p}^{2}}, \v x_{k} - \v x_{\av}>\\
&&\displaystyle  +  \vn{\frac{\partial(\v x_{k}-\v x_{\av})}{\partial \theta_{p}}}^{2} \\
&&\displaystyle +   < \frac{\partial^{2}(\v x_{k}-\v x_{\av})}{\partial \theta_{p}^{2}}, (\v r_{k} - \v r_{\av})\times \v \omega_0>. 
\end{array}
\]
Let $\bar{\v q}= (\bar{\v x}, \bar{\v r})$ be a point belonging to $\Lambda$. 
By using the relations (\ref{condpro}) (characterizing the set $\Lambda$) we observe that in the set $\Lambda$ the following conditions hold 
\[
\begin{array}{rcll }
\v x_{k} - \v x_{\av} &=& (\v r_{k}-\v r_{\av})\times \v \omega_0, &\quad \quad k \in G_{1}\\
(\v r_{k}-\v r_{\av})\times \v \omega_0 &=&0,&\quad \quad  k \in G_{2}.\\
\end{array}
\]
this yields
\begin{eqnarray}
\displaystyle \left. \frac{\partial^{2}S}{\partial \theta_{p}^{2}} \right |_{\bar{\v q}} 
& =& \displaystyle \frac{N-1}{N^{2}}\vn{\frac{\partial \v x_{p}}{\partial \theta_{p}}}^{2} + \frac{N-1}{N}< \frac{\partial^{2}\v x_{p}}{\partial \theta^{2}_{p}},\v x_{p}-\v x_{\av}>  \nonumber\\&&+ \frac{(N-1)^{2}}{N^{2}}\vn{\frac{\partial \v x_{p}}{\partial \theta_{p}}}^{2} \nonumber
\\&&\displaystyle - \sum_{k\in G_{2}\setminus p} \frac{1}{N}<\frac{\partial^{2} \v x_{p}}{\partial \theta_{p}^{2}}, \v x_{k}-\v x_{\av}>. \label{eq:S_dd}
\end{eqnarray}
Since $\displaystyle \v x_{av} = \alpha \frac{\v \omega_0}{\vn{\v \omega_0}},\, 0 \leq \alpha <1$ and $\displaystyle \frac{\partial^{2} \v x_{p}}{\partial \theta_{p}^{2}} = \displaystyle \frac{\v \omega_0}{\vn{\v \omega_0}}$ in $\Lambda$, the expression (\ref{eq:S_dd}) reduces to
\[
\begin{array}{rcl}
\displaystyle
\left. \frac{\partial^{2}S}{\partial \theta_{p}^{2}} \right |_{\bar{\v q}} &=&\displaystyle \frac{N-1}{N^{2}} -(\alpha +1) \frac{N-1}{N} + \frac{(N-1)^{2}}{N^{2}}\\
&&\displaystyle - \alpha\frac{(N-M-1)}{N} - \frac{1}{N}<0,
\end{array}
\]  
which shows that (\ref{S}) does not attain a local minimum in the set $\Lambda$. Let $\Lambda_{\bar{\v q}}$ be the connected component of $\Lambda$ containing $\bar{\v q}$. Consider a neighborhood $B(\bar{\v q})$ in the shape space such that $B(\bar{\v q})\setminus \Lambda_{\bar{\v q}}$ contains no points where $\dot S =0$. Choose a point $\tilde{q} \in B(\bar{\v q})$ such that $S(\tilde{\v q})<S(\bar{\v q})$. Since the function $S$ decreases along the solutions, the solution with initial condition $\tilde{\v q}$ cannot converge to $\Lambda_{\bar{\v q}}$ and leaves  
$B(\bar{\v q})$ after a finite time. Since $S$ is not at a local minimum in $\Lambda_{\bar{\v q}}$ we can take $\tilde{\v q}$ arbitrary close to $\bar{\v q}$ which shows that $\bar{\v q}$ is unstable.

We conclude that the set $E(\v \omega_0)$ is asymptotically stable in the shape space and that the other positive limit sets are unstable. 
\qedp

\section{Proof of Theorem 3}
The function $B\triangleq Q+S$ is non negative and it is non-increasing along the
solutions of (\ref{model}) with the control law (\ref{concon4}).
Then $B$ converges to a limit as $t \rightarrow \infty$.
Furthermore the second derivative $\ddot B$ is bounded (because $\v v_k^{a} -  \v v^a_{\av}$ is bounded for every $k$). From Barbalat's Lemma $\dot B \rightarrow 0$ when $t\rightarrow \infty$ and therefore the solutions converge to the set where  
\begin{equation}\label{bsymm1}
\left[\v v^{a}_{k} -\v v^a_{\av} + \left(\frac{<\v \omega_0, \v x_{\av}>}{\vn{\v \omega_0}}-\alpha\right)\frac{\v \omega_0}{\vn{\v \omega_0}} \right]\times \v x_k=0.
\end{equation}
The $\v x$ dynamics in this set reduce to
\[
\dot{\v x}_{k} =\v  \omega_0 \times \v x_{k} \quad \quad k=1,\ldots, N.
\]
Following the same lines of the proof of Theorem $2$, we analyze the stability of the positive limit sets. \\ 
i) Suppose that $\v \omega_0 \times \v x_{k} \neq 0$ for every $k$. Then the only possible way for (\ref{bsymm1}) to hold is that $\v v_k^a -\v v^a_{\av} + (\frac{<\v \omega_0, \v x_{\av}>}{\vn{\v \omega_0}}-\alpha)\frac{\v \omega_0}{\vn{\v \omega_0}} =0$ for every $k$. Factoring the first term in parallel and orthogonal components (with respect to $\v \omega_0$) we obtain
\[
\begin{array}{c}
<\v x_{k}-\v x_{\av},\v \omega_0>\frac{\v \omega_0}{\vn{\v \omega_0}^{2}} + \frac{\v \omega_0}{\vn{\v \omega_0}^{2}} \times ((\v v_k^a-\v v^a_{\av}) \times \v \omega_0)\\ 
+ \left(\frac{<\v \omega_0, \v x_{\av}>}{\vn{\v \omega_0}}-\alpha\right)\frac{\v \omega_0}{\vn{\v \omega_0}} =0,
\end{array}
\]
which implies that $<\v x_{k},\v \omega_0>\frac{\v \omega_0}{\vn{\v \omega_0}}= \alpha$ and 
$ \v v_k^a=\v v^a_{\av}$ for every $k$. The second condition tell us that a relative equilibrium is reached while the first says that the pitch of every particle is fixed to the desired value $\alpha$. Since in this set the Lyapunov function attains a global minimum we conclude that the set of relative equilibria with rotation vector $\v \omega_0$ and pitch $\alpha$ is asymptotically stable in the shape space.\\ 
ii) Suppose that $\v \omega_0 \times \v x_k \neq 0$ for $k \in G_1$ and $\v \omega_0 \times \v x_j =0$
for $j \in G_2$, where $G_{1}$ and $G_{2}$ are defined in the proof of Theorem $2$. 
In such a configuration we obtain
\begin{equation*}%\label{condpro}
\begin{array}{rcl}
\v v_k^a - \v v^a_{\av} + (\frac{<\v \omega_0, \v x_{\av}>}{\vn{\v \omega_0}}-\alpha)\frac{\v \omega_0}{\vn{\v \omega_0}} &=& 0, \quad \quad k\in G_1\\
(\v v^a_j - \v v^a_{\av})\times \v \omega_0  &=& 0,  \quad \quad j\in G_2,
\end{array}
\end{equation*}
where $\v v^{a}_j\neq \v v^a_{\av},\, j\in G_2$.
 Following the same lines of the proof of Theorem 2, it can be shown (by calculating the second derivative of the Lyapunov function with respect to a suitable direction) that the set defined by this configuration is unstable (unless $|G_{2}| = 0$ that is the case considered in the point i)).         \qedp
\section{Proof of Theorem 4}
To show that the resulting closed-loop
vector field is invariant under an action of $SE(3)$, it is
sufficient to observe that the dynamic control law (\ref{concon5})
depends only on the relative orientations and relative positions
of the particles. With the change of variables $\v \omega_k^a =
R_k \v \omega_k$ (\ref{concon5}) rewrites to
\begin{subequations}\label{inv}
\begin{eqnarray}
\v u_k &=&  R_k^T \left(\v \omega_k^a + \left[\left(\v r_k - \v r_{\av}\right)\times \v \omega_k^a -\v x_{\av}\right] \times \v x_k\right)\hspace{-1 mm},\label{inv1}\\
 \dot{\v \omega}_k^a &=& \sum_{j=1}^N \left(\v \omega^a_j
- \v \omega^a_k\right). \label{inv2}
\end{eqnarray}
\end{subequations}
We observe that (\ref{inv2}) is independent of the particle
dynamics. Therefore the solutions of (\ref{inv2}) will exponentially converge to
a consensus value ${\v \omega}_{\av}\triangleq\frac{1}{N}\sum_{j=1}^{N} R_{j}(0)\v \omega_{j}(0)$, i.e. ${\v \omega}_k^a
\rightarrow {\v \omega}_{\av}$ when $t\rightarrow \infty$, for every
$k=1,2,\ldots,N$. 
Therefore (\ref{inv1})
asymptotically converge to
\begin{equation}\label{asymp}
\v u_k = R_k^T \left({\v \omega}_{\av} + \left[\left(\v r_k - \v r_{\av}\right)\times {\v \omega}_{\av} -\v x_{\av}\right] \times \v
x_k\right).
\end{equation}
The positive limit sets (in the shape space) for system (\ref{model}) with the control law (\ref{asymp}) have been  analyzed in Theorem 2 and we already know that $E(\v \omega_{\av})$ is an asymptotically stable set. Therefore system (\ref{model}) with (\ref{inv}) is a cascade of an exponentially stable system and a system with an asymptotically stable set (in the shape space) $E(\v \omega_{\av})$. From standard results (see e.g. \cite{SeJaKo,So89}) we conclude that $E \times C_{\v \omega}$ is a stable attractor, in the shape space, for the cascade system. The instability of the other positive limit sets follows from Theorem 2.    
\qedp

\section{Proof of Proposition 2}
Following the same lines of the proof of Theorem 2 we observe that the only asymptotically stable equilibria of the $\v v$ dynamics are relative equilibria of  (\ref{model}). These configurations are characterized by $\v v_k^a = \tilde{\v v}_{\av}, k=1,2,\ldots, N$. Since  $ \tilde{\v v}_{\av} = \frac{1}{N+1}\left(\sum_{k=1}^{N} \v v_k^a +\v v_{0}\right)$, we conclude that  $\v v_k^a = \v v_{0},\, k=1,2,\ldots, N$.  \qedp
 
\section{Proof of Theorem 6}
Since the control law does not depend on the relative spacing, we analyze the reduced dynamics on relative orientations.
Set $\v b^{a}_k=R_k \v b_k $. Then $\v b^{a}(t)$ obeys the
consensus dynamics $\dot{\v b}^a=-\tilde L(t)\v b^{a}$, which implies
that its solutions exponentially converge to a consensus value $\v b_{0}$.
Therefore, the control law
\begin{equation}\label{as_aut}
\v u_k=R_k^T(\v x_k \times \v b_k^{a}),
\end{equation}
asymptotically converges to the control
\begin{equation}\label{aut}
\v u_k=R_k^T(\v x_k \times \v b_0),
\end{equation}
for every $k=1,2,\ldots,N$.
The limiting system is decoupled
into $N$ identical systems whose limit sets (of the reduced dynamics) are
characterized by $\v x_k = \frac{\v b_0}{\vn{\v
b_0}}$ or $\v x_k = -\frac{\v b_0}{\vn{\v b_0}}$
for every $k$. The synchronized set $\v x_k = \frac{\v
b_0}{\vn{\v b_0}}$ is exponentially stable while the set characterized by $\v x_k = -\frac{\v b_0}{\vn{\v b_0}}$ is
unstable. Therefore system (\ref{model}) with (\ref{inv}) is a cascade of a uniformly exponentially stable system (in the shape space) with a system with an asymptotically stable set (in the shape space). From standard results on stability of cascade systems, we conclude that $E(0)\times C_{\v b}$ is a stable attractor, in the shape space, for the cascade system. The instability of the other positive limit sets follows from the instability of the corresponding limit sets in the (limit) decoupled dynamics.   \qedp

\section{Proof of Theorem 7}
Observe that with the change of variables $\v \omega_k^a=R_k \v \omega_k,\, \v b_k^a=R_k \v b_k,\, \v
c_k^a=R_k \v c_k + \v r_k$ (\ref{concon6}) rewrites to
\[
\begin{array}{rcl}
\v u_k &=&  R_{K}^{T}\left(\omega^{a}_k + \left[(\v r_k - \v c_k^a)\times \v \omega_k^a -\v b_k^a \right] \times \v x_k\right)\\
\dot{\v \omega}_k^a &=&  - \sum_{j=1}^N
L_{kj}(t) \v \omega^a_j\\
\dot{\v b}_k^a &=&- \sum_{j=1}^N L_{kj}(t)
\v b_j^a\\
\dot{\v c}_k^a &=& -\sum_{j=1}^N  L_{kj}(t) \v c_j^a
\end{array}
\]
and the consensus dynamics are not influenced by the particles dynamics. Therefore, from Theorem (\ref{ConsEucl}), 
%\ref{ConsEucl}
we conclude that the variables $\v \omega_k^a,\v b_k^a$ and $\v c_k^a$ asymptotically converge to the consensus values ${\v
\omega}_{0},{\v b}_{0}$ and ${\v c}_{0}$ respectively, and the particles' dynamics become asymptotically decoupled.  %System (\ref{model}) equipped with the control law (\ref{concon6}) will
%converge to an asymptotically autonomous system whose limit set is characterized in Theorem \ref{th:1}.  
The dynamics of the decoupled system can be easily characterized defining the Lyapunov function
\[
\tilde V = \sum_{k=1}^N \vn{\v v_k^a - \v v_0}^2, 
\] 
where $\v v_0 = \v c_0\times \v \omega_0 + \v b_0$ and $\v v_k^a = \v r_k\times \v \omega_0 + \v x_k$. Now observe that $\tilde V$ is non increasing along the solutions of the decoupled system:
\[
\dot {\tilde V} = -\sum_{k=1}^N \vn{(\v v_k^a - \v v_0)\times \v x_k}^2,
 \]
that is sufficient to conclude that the set of relative equilibria with rotation vector $\v \omega_0$ and $\v v_k^a =\v  v_0,\, k=1,2,\ldots,N$  is asymptotically stable for the uncoupled dynamics.   
Following the same lines of the proof of Theorem 6, we conclude that the set $E \times C_{\v \omega}\times C_{\v b}\times C_{\v c}$ is asymptotically stable in the shape space. The instability of the other limit sets follows from the instability of the corresponding limit sets in the (limit) decoupled dynamics. \qedp


\begin{thebibliography}{10}
\providecommand{\url}[1]{#1}
\csname url@samestyle\endcsname
\providecommand{\newblock}{\relax}
\providecommand{\bibinfo}[2]{#2}
\providecommand{\BIBentrySTDinterwordspacing}{\spaceskip=0pt\relax}
\providecommand{\BIBentryALTinterwordstretchfactor}{4}
\providecommand{\BIBentryALTinterwordspacing}{\spaceskip=\fontdimen2\font plus
\BIBentryALTinterwordstretchfactor\fontdimen3\font minus
  \fontdimen4\font\relax}
\providecommand{\BIBforeignlanguage}[2]{{%
\expandafter\ifx\csname l@#1\endcsname\relax
\typeout{** WARNING: IEEEtran.bst: No hyphenation pattern has been}%
\typeout{** loaded for the language `#1'. Using the pattern for}%
\typeout{** the default language instead.}%
\else
\language=\csname l@#1\endcsname
\fi
#2}}
\providecommand{\BIBdecl}{\relax}
\BIBdecl

\bibitem{LePaLeFrSeDa}
N.~Leonard, D.~Paley, F.~Lekien, R.~Sepulchre, D.~Fratantoni, and R.~Davis,
  ``Collective motion, sensor networks and ocean sampling,'' \emph{Proceedings
  of the IEEE}, vol.~95, no.~1, pp. 48--74, 2007.

\bibitem{bullo}
J.~Cortes, S.~Martinez, T.~Karatas, and F.~Bullo, ``Coverage control for mobile
  sensing networks,'' \emph{IEEE Trans. on Robotics and Automation}, vol.~20,
  no.~2, pp. 243--255, 2004.

\bibitem{JuKr_cdc}
E.~W. Justh and P.~S. Krishnaprasad, ``Natural frames and interacting particles
  in three dimensions,'' in \emph{Proceedings of the 44th IEEE Conference on
  Decision and Control and European Control Conference}, Seville, Spain, 2005,
  pp. 2841--2846.

\bibitem{FaMu}
J.~Fax and R.~Murray, ``Information flow and cooperative control of vehicle
  formations,'' \emph{IEEE Trans. on Automatic Control}, vol.~49, no.~9, pp.
  1465--1476, 2004.

\bibitem{CoKrJaRuFr}
I.~Couzin, J.~Krause, R.~James, G.~Ruxton, and N.~Franks, ``Collective memory
  and spatial sorting in animal groups,'' \emph{J. Theor. Biol.}, vol. 218, pp.
  1--11, 2002.

\bibitem{Brockett76}
R.~W. Brockett, ``Nonlinear systems and differential geometry,''
  \emph{Proceedings of the IEEE}, vol.~64, pp. 61--72, 1976.

\bibitem{Brockett72}
------, ``System theory on group manifolds and coset spaces,'' \emph{SIAM
  Journal on Control}, vol.~10, pp. 265--284, 1972.

\bibitem{Brockett73b}
------, \emph{Lie Algebras and Lie Groups in Control Theory}.\hskip 1em plus
  0.5em minus 0.4em\relax Dordrecht, The Netherlands: Reidel Publishing Co.,
  1973, pp. 43--82.

\bibitem{Brockett73a}
------, ``Lie theory and control systems defined on spheres,'' \emph{SIAM
  Journal Appl. Math.}, vol.~25, pp. 213--115, 1973.

\bibitem{Wood}
R.~W. Brockett and J.~R. Wood, \emph{Electrical Networks Containing Controlled
  Switches}.\hskip 1em plus 0.5em minus 0.4em\relax Western Periodicals Co.,
  1974, pp. 1--11.

\bibitem{BlochBook}
A.~M. Bloch, \emph{Nonholonomic Mechanics and Control}.\hskip 1em plus 0.5em
  minus 0.4em\relax New York: Springer-Verlag, 2003.

\bibitem{Navin}
N.~Khaneja, ``Switched control of electron nuclear spin systems,'' \emph{Phys.
  Rev. A}, vol.~76, p. 012316, 2007.

\bibitem{ScSe_cdc06}
L.~Scardovi and R.~Sepulchre, ``Collective optimization over average
  quantities,'' in \emph{Proceedings of the 45th IEEE Conference on Decision
  and Control}, San Diego, Ca, 2006, pp. 3369--3374.

\bibitem{scardovi:2007}
L.~Scardovi, A.~Sarlette, and R.~Sepulchre, ``Synchronization and balancing on
  the {$N$}-torus,'' \emph{Systems and Control Letters}, vol.~56, no.~5, pp.
  335--341, 2007.

\bibitem{FrYaLy}
R.~A. Freeman, P.~Yang, and K.~M. Lynch, ``Distributed estimation and control
  of swarm formation statistics,'' in \emph{Proceedings of the American Control
  Conference}, Minneapolis, MN, 2006, pp. 749--755.

\bibitem{Ol}
R.~Olfati-Saber, ``Distributed kalman filter with embedded consensus filters,''
  in \emph{Proceedings of the 44th IEEE Conference on Decision and Control and
  European Control Conference}, Seville, Spain, 2005, pp. 8179--8184.

\bibitem{SePaLe_limited}
R.~Sepulchre, D.~Paley, and N.~Leonard, ``Stabilization of planar collective
  motion with limited communication,'' \emph{IEEE Trans. on Automatic Control},
  vol.~53, no.~3, pp. 706--719, 2008.

\bibitem{SePaLe}
------, ``Stabilization of planar collective motion: All-to-all
  communication,'' \emph{IEEE Trans. on Automatic Control}, vol.~52, no.~5, pp.
  811--824, 2007.

\bibitem{ecc}
L.~Scardovi, R.~Sepulchre, and N.~Leonard, ``Stabilization laws for collective
  motion in three dimensions,'' in \emph{European Control Conference}, Kos,
  Greece, 2007, pp. 4591--4597.

\bibitem{ScLeSe}
L.~Scardovi, N.~Leonard, and R.~Sepulchre, ``Stabilization of collective motion
  in three dimensions: A consensus approach,'' in \emph{Proceedings of the 46th
  IEEE Conference on Decision and Control}, New Orleans, La, 2007, pp.
  2931--2936.

\bibitem{SaSeLe}
A.~Sarlette, R.~Sepulchre, and N.~Leonard, ``Autonomous rigid body attitude
  synchronization,'' in \emph{Proceedings of the 46th IEEE Conference on
  Decision and Control}, New Orleans, La, 2007, pp. 2566--2571.

\bibitem{IgFuSp}
Y.~Igarashi, M.~Fujita, and M.~W. Spong, ``Passivity-based 3d attitude
  coordination: Convergence and connectivity,'' in \emph{Proceedings of the
  46th IEEE Conference on Decision and Control}, New Orleans, La, 2007, pp.
  2558--2565.

\bibitem{MuZeSh}
R.~Murray, Z.~Li, and S.~S. Sastry, \emph{A Mathematical Introduction To
  Robotic Manipulation}.\hskip 1em plus 0.5em minus 0.4em\relax CRC Press,
  1994.

\bibitem{Mocdc}
L.~Moreau, ``Stability of continuous-time distributed consensus algorithms,''
  in \emph{Proceedings of the 43rd IEEE Conference on Decision and Control},
  Paradise Island, Bahamas, 2004, pp. 3998--4003.

\bibitem{Mo}
------, ``Stability of multi-agent systems with time-dependent communication
  links,'' \emph{IEEE Trans. on Automatic Control}, vol.~50, no.~2, pp.
  169--182, 2005.

\bibitem{OlMu}
R.~Olfati-Saber and R.~Murray, ``Consensus problems in networks of agents with
  switching topology and time-delays,'' \emph{IEEE Trans. on Automatic
  Control}, vol.~49, no.~9, pp. 1520--1533, 2004.

\bibitem{phdbhatta}
P.~Bhatta, ``Nonlinear stability and control of gliding vehicles,'' Ph.D.
  dissertation, Princeton University, 2006.

\bibitem{SeJaKo}
R.~Sepulchre, M.~Jankovi\'c, and P.~Kototovi\'c, \emph{Costructive Nonlinear
  Control}.\hskip 1em plus 0.5em minus 0.4em\relax Springer, 1997.

\bibitem{So89}
E.~D. Sontag, ``Remarks on stabilization and input-to-state stability,'' in
  \emph{Proceedings of the 46th IEEE Conference on Decision and Control},
  Tampa, FL, 1989, pp. 1376--1378.

\end{thebibliography}
\end{document}